\documentclass[10pt]{book}
\usepackage{array}
\usepackage{amsfonts}
\begin{document}
\renewcommand{\baselinestretch}{1.2}
\markright{
}
\markboth{\hfill{\footnotesize\rm DOROTA M. DABROWSKA AND WAI TUNG HO
}\hfill}
{\hfill {\footnotesize\rm A MODULATED RENEWAL PROCESS} \hfill}
\renewcommand{\thefootnote}{}
$\ $\par
\fontsize{10.95}{14pt plus.8pt minus .6pt}\selectfont
\vspace{0.8pc}
\centerline{\large\bf ESTIMATION IN A SEMIPARAMETRIC }
\vspace{2pt}
\centerline{\large\bf MODULATED RENEWAL PROCESS}
\vspace{.4cm}
\centerline{Dorota M. Dabrowska}
\vspace{.4cm}
\centerline{\it University of California, Los Angeles} 
\vspace{.4cm}
\centerline{Wai Tung Ho}
\vspace{.4cm}
\centerline{\it SPSS Inc., Chicago} 
\vspace{.55cm}

\fontsize{9}{11.5pt plus.8pt minus .6pt}\selectfont

\begin{quotation}
\noindent {\it Abstract:}

We consider parameter estimation 
in a regression model corresponding to  an iid sequence 
of censored observations of a  
finite state modulated renewal process. The model assumes a 
similar form as in  Cox regression except that the baseline
intensities  are functions of the backwards recurrence time
of the
process and a time dependent covariate. As a result of this
it falls outside the class of multiplicative intensity 
counting process models. We use kernel estimation to construct
estimates of the regression coefficients and baseline cumulative 
hazards. We give conditions for 
consistency and asymptotic normality of estimates. Data from 
a bone marrow transplant study are used to illustrate the results. 

\vspace{9pt}
\par
\end{quotation}\par

\fontsize{10.95}{14pt plus.8pt minus .6pt}\selectfont

\voffset = .5in

\renewcommand{\thesection}{\arabic{section}}
\renewcommand{\theequation}{\arabic{section}.\arabic{equation}}

\renewcommand{\theenumi}{\arabic{enumi}}
\renewcommand{\theenumii}{\arabic{enumii}}

\font\bbigsym=cmmi10 scaled\magstep4 
\font\bbbigsym=cmmi10 scaled\magstep5 
\newcommand{\BIGP}{\lower2pt\hbox{\bbigsym\char'031}} 
\newcommand{\BBIGP}{\lower3pt\hbox{\bbbigsym\char'031}} 
 
\newcommand{\hash}{{\rm\char'043}} 
\newcommand{\openE}{\hbox{I\kern-.2em E}} 
\newcommand{\openN}{\hbox{I\kern-.2em N}} 
\newcommand\openR{\hbox{I\kern-.2em R}} 
\newcommand\smopenR{{\hbox{\sevenrm I \kern-.5em R}}} 
\newcommand\openH{\hbox{I\kern-.2em H}}

\newcommand{\empE}{\mathbb E\, }
\newcommand{\empP}{\mathbb P}
\newcommand{\empG}{\mathbb G}
\newcommand{\empU}{\mathbb U}
\newcommand{\empJ}{\mathbb J}
 
\def\iint{\mathop{\hbox{$\displaystyle\int\int$}}} 
\def\idotsint{\mathop{\hbox{$\displaystyle\int\dots\int$}}} 
\def\convP{\buildrel\P\over\to} 
\def\convD{\buildrel\Cal D\over\to}

\renewcommand\epsilon{\varepsilon} 
\renewcommand \phi{\varphi}
\newcommand \ow{\omega}
\newcommand \Ow{\Omega}
\newcommand \eps{\epsilon}
 
\renewcommand{\hat}{\widehat} 
\renewcommand\tilde{\widetilde} 
\renewcommand{\bar}{\overline}

\newcommand\implies{\Rightarrow} 
\newcommand\Verti{\Vert_{\infty}}
\newcommand\Vertv{\Vert_v}
\newcommand\Vertd{\Vert_{\delta}}
\newcommand\Vertt{\Vert_{t_0}} 
\newcommand\Pern{\text{Per}_n}
\renewcommand{\Pr}{{\rm Pr}}
\newcommand\Bor{\text{Bor}}
\newcommand\Exp{\text{Exp}}
\newcommand\EXP{\text{EXP}}
\newcommand\Peano{{\cal P}}
\newcommand{\corr}{\rm corr}
\newcommand{\cov}{{\rm cov}}
\newcommand{\var}{{\rm Var} }
\newcommand{\Var}{{\rm Var} }
\newcommand{\sign}{{\rm sign}}
\newcommand{\Vec}{{\rm Vec}}
\newcommand{\diag}{{\rm diag}}
\newcommand{\done}{$ \ \Box$ }

\newcommand{\A}{{\cal A}}
\newcommand{\B}{{\cal B}}
\newcommand{\C}{{\cal C}}
\newcommand{\D}{{\cal D}}
\newcommand{\E}{{\cal E}}
\newcommand{\F}{{\cal F}}
\newcommand{\G}{{\cal G}}
\renewcommand{\H}{{\cal H}}
\newcommand{\I}{{\cal I}}
\newcommand{\J}{{\cal J}}
\newcommand{\K}{{\cal K}}
\renewcommand{\L}{{\cal L}}
\newcommand{\M}{{\cal M}}
\newcommand{\N}{{\cal N}}
\renewcommand{\P}{{\cal P}}
\newcommand{\R}{{\cal R}}
\renewcommand{\S}{{\cal S}}
\newcommand{\T}{{\cal T}}
\newcommand{\U}{{\cal U}}
\newcommand{\V}{{\cal V}}
\newcommand{\W}{{\cal W}}
\newcommand{\X}{{\cal X}}
\newcommand{\Y}{{\cal Y}}
\newcommand{\Z}{{\cal Z}}

\newcommand{\al}{{\alpha}}
\newcommand{\lam}{{\lambda}}
\renewcommand{\l}{{\ell}}
 
\newcommand\cldot{\cdot \ldots \cdot}  

\def\for{\quad {\rm for} \quad} 
\def\and{\quad {\rm and} \quad}
\def\tor{\quad {\rm or} \quad}
\def\tif{\quad {\rm if} \quad}
\def\tiff{\quad{\rm iff} \quad}

\begingroup
\newtheorem{bigthm}{Theorem}
\newtheorem{thm}{Theorem}[section]
\newtheorem{prop}[thm]{Proposition}
\newtheorem{lem}[thm]{Lemma}
\newtheorem{cor}[thm]{Corrolary}
\newtheorem{defn}[thm]{Definition}
\newtheorem{ex}[thm]{Example} 
\newtheorem{rem}[thm]{Remark} 
\newtheorem{con}[thm]{Condition} 
\endgroup

\newcommand\bthm {\begin{thm} \rm}
\newcommand\ethm { \end{thm}}
\newcommand\bprop {\begin{prop} \rm}
\newcommand\eprop { \end{prop}}
\newcommand\blem {\begin{lem} \rm}
\newcommand\elem { \end{lem}}
\newcommand\bcor {\begin{cor} \rm}
\newcommand\ecor { \end{cor}}
\newcommand\bdf {\begin{defn} \rm}
\newcommand\edf { \end{defn}}
\newcommand\bex {\begin{ex} \rm}
\newcommand\eex { \end{ex}}
\newcommand\brem {\begin{rem} \rm}
\newcommand\erem { \end{rem}}
\newcommand\bcon {\begin{con} \rm}
\newcommand\econ { \end{con}}

\newcommand\bitem { \begin{itemize}}
\newcommand\eitem { \end{itemize}}

\newcommand\beq {\begin{equation}}
\newcommand\eeq {\end{equation}}

\newcommand\bea{\begin{eqnarray}}
\newcommand\eea{\end{eqnarray}}

\newcommand\beaa{\begin{eqnarray*}}
\newcommand\eeaa{\end{eqnarray*}}

\newcommand{\Proof} {\noindent \it Proof . \rm}

\newcommand{\etal} {\it et al. \rm}

\newcommand\bib{\bibliography}

\newcommand\bbib{\begin{thebibliography}{}}
\newcommand\ebib{\end{thebibliography}{}}
\newcommand\bibit{\bibitem{}}

\newcommand\bref{\section*{References} 
\begin{list}{0}{\setlength{\rightmargin}{\leftmargin}}}

\newcommand\eref{\end{list}}

\newcommand{\sect}[1]{\section{#1}\setcounter{equation}{0}}
\newcommand{\subsect}[1]{\subsection{#1}}
\newcommand{\subsubsect}[1]{\subsubsection{#1}}

\newpage

\setcounter{section}{1}
\setcounter{equation}{0} 
\noindent {\bf 1. Introduction}

In medical and engineering applications it is common to consider
a Markov renewal process to model the lengths of time spent
in consecutive stages of a disease or lifetime of a
piece of equipment. Denoting by $\empJ = \{1, \ldots,k\}$ the set
of possible states, the process is described by a sequence of
random variables $(T, J) = (T_m,J_m)_{m \geq 0}$, such that 
$T_0 < T_1 < T_2 < \ldots$ are consecutive times of entrances
into states $J_0, J_1, \ldots, J_m \in \empJ$. Under assumption of the
Markov renewal process,  the sequence $J = \{J_m: m \geq 0\}$ of
states visited forms a Markov chain and given $J$,
the sojourn times $T_1, T_2 - T_1, \ldots$ are independent
with distributions depending on the adjoining states. 
Associated with the sequence $(T, J)$ is  a 
counting process $\{\tilde N_{ij}(t): t \geq 0, i,j \in \empJ \}$
whose components register each direct $i \to j$ transition,
$$  
\tilde N_{ij}(t) = 
\sum_{m \geq 0}1(T_{m+1} \leq t, J_{m+1} = i, J_m = j) \;.
$$
Its compensator $\{\Lambda_{ij}(t): t \geq 0, i,j \in \empJ\}$
relative to the self--exciting filtration is given by
$$
\Lambda_{ij}(t) = \int_0^t 1(J(s-) = i) dA_{ij}(L(s)) \;,
$$
where $J(t), t \geq 0$ is the state occupied at time t, 
$L(t) = t - T_{\tilde N(t-)}, \tilde N(t) = \sum_{ij \in E } 
\tilde N_{ij}(t)$ is the
backwards recurrence time, and $[A_{ij}(x)]_{i,j \in \empJ}$ is a matrix of
unknown deterministic functions representing cumulative hazards of one-step
transitions. Nonparametric estimation of this matrix and the associated
semi-Markov kernel of the process was considered by 
Lagakos, Sommer and Zelen (1978), Gill (1980), Voelkel and Crowley (1984),
and Phelan (1990), among others.

In this paper we consider estimation in a modulated renewal 
process, assuming that  components of the counting process
$\{\tilde N_{ij}: (i,j) \in \empJ \}$ have intensities of the form
\beq
\Lambda_{ij}(t) = \int_0^t 1(J(s-)=i) e^{\beta^T Z_{ij}(s)} 
\alpha_{ij}(L(s), X(s))
ds \;,
\eeq
where 
$X(s)$ is a time dependent covariate,
$Z = \{Z_{ij}(t): t \geq 0, i,j \in \empJ\}$ is a vector of external
transition specific covariates, and $[\alpha_{ij}]$ is a matrix
of two-parameter baseline hazards.
 A model of this kind may
arise for instance in medical applications where survival status 
of a patient is characterized by an illness process with baseline
intensities dependent on the length of time spent in each stage of a 
disease and a covariate $X(s)$, possibly changing with time. 
In the absence of this covariate, the model reduces to the modulated
renewal process proposed by Cox (1973) with cumulative intensities
\beq
\Lambda_{ij}(t) = \int_0^t 1(J(s-) = i) e^{\beta^T Z_{ij}(s)} 
\alpha_{i,j}(L(s)) ds \;.
\eeq
Both models have
several interesting features. The
first one is that the 
event times can be viewed as recorded on two simultaneously
evolving time scales. In the case of  (1.2), 
the covariates  depend on the calendar time $t$,
whereas the matrix $\alpha$ of baseline hazards depends on the duration scale. 
In the case of  (1.1), the latter matrix depends both on the duration
and calendar time scale.
Further, if $\alpha$ corresponds to a matrix of functions depending 
only on a Euclidean parameter $\theta$, then estimation of the pair 
$(\beta, \theta)$
based on an iid sample of modulated renewal processes 
can be carried out using a counting process
framework for analysis of maximum likelihood or M estimates. However, if the
matrix $\alpha$ is completely unspecified, 
then its nonparametric maximum likelihood 
estimate falls outside the class of statistics taking the form 
of stochastic integrals
with respect to counting processes (Gill (1980)). Similarly,
in the case of  (1.2), estimation of the regression 
coefficient $\beta$ can be in principle based on the solution to the score
equation 
\beq
\Phi_n(\beta) = \sum_{\ell=1}^n \sum_{i,j \in \empJ} \int 
[ Z_{ij\l}(t) - {S^{(1)} \over S^{(0)}}(t, \beta) ] 
\tilde N_{ij\l}(dt) = 0 \;,
\eeq
where $S^{(p)}(t, \beta) =
\sum_{\ell= 1}^n 1(J_{\ell}(t-) = i) Z_{ij\ell}^{p}(t)e^{\beta Z_{ij\ell}(t)}, p = 0,1$. 
However, as a result of the
dependence of the compensators on the backwards recurrence time, 
the  score function  in (1.3),
evaluated at the true parameter value $\beta_0$, fails to satisfy the identity
$\empE \Phi_n(\beta_0) = o_P(1)$, and consequently the estimate of the 
regression coefficient obtained by solving the equation $\Phi_n(\beta) = 0$
cannot be consistent.  Several authors considered also the special case of
the one-jump process (1.1)  and showed that estimation 
of regression coefficients requires smoothing (Sasieni (1992), 
Dabrowska, (1997),
Nielsen, Linton and Bickel (1998), Pons and Vissier, (2000)).

To circumvent difficulties arising in the analyses of renewal
processes,  Gill (1980) and Oakes and Cui (1994)
proposed the use of a random 
time-change approach which replaces the calendar time scale $t$ by 
the duration scale. 
Here we consider an extension of this approach to 
analyse  a simple case of  (1.1),  assuming that the 
covariate 
$X(s)$ is constant
between the jumps of the process $\tilde N(t) = \sum_{ij} \tilde N_{ij}(t)$,
and $\{Z_{ij}(t): i,j \in \empJ, t \geq 0\}$ is a vector of
external covariates. In Section 3 we discuss kernel estimation 
in single-type models. In Section 4 we give examples  multi-type
models with a ``small'' state space to which the results can also
be applied. We use data from a bone marrow transplant
study to illustrate the results.

\par

\setcounter{section}{2}
\setcounter{equation}{0} 
\noindent{\bf 2. The model}

Throughout the paper we assume that $(\Omega, \F, P)$ is a complete
probability space and $(T_m, V_m)_{m \geq 0}$ is a marked point process
defined on it with marks taking on values in a measurable space $(E, \E)$
and enlarged by the empty mark $\Delta$.
Thus $T_0 < T_1 < \ldots T_m \ldots $ is a sequence of random
time points registering occurrence of some events in time, and
such that $T_m$ are almost surely distinct and $T_m \uparrow \infty$
P-a.s. At time $T_m$ we observe a variable $V_m$ such that
$V_m  \in  E$ if $T_m < \infty$, and $V_m  =  \Delta$ if $T_m = \infty$.

For any $B \in \E$, let 
$\tilde N(t,B) = \sum_{m \geq 0} 1(T_{m+1} \leq t, V_{m+1} \in B)$
be the process counting observations falling into the set $[0, t] \times B$.
The internal history of the process, $\{\F_t^N \}_{ t \geq 0}$, 
represents information
collected on $N$ until time $t$, and is given by
$$
\F_t^N = \sigma(1(T_m \leq s, V_m \in B): m \geq 0, s \leq t, B \in \E) \;.
$$
Then $\{\F_t^N \}_{t \geq 0}$ forms an increasing  family of 
right-continuous $\sigma$-fields.
Let $\F_t = \F_0 \vee \F_t^N$ be the self-exciting filtration associated
with the process $\tilde N$, 
obtained by adjoining to the internal history of the
process, the $P$-null sets. The compensator of 
the process $\tilde N(t,B)$, with respect
to $\F_t$ is given by
$$
\Lambda(t,B) = \Lambda(T_m, B) + \int_{(T_m, t]} {P_m(d(s,v)) \over 
P_m([s, \infty);E \cup \Delta)}  \for t \in (T_m, T_{m+1}] \;,
$$
where $P_m(d(s,v))$ is a version of a regular conditional distribution
of $(T_{m+1},V_{m+1})$ given $\F_{T_m}$ (Jacod (1975)).

In this paper we  assume that the marks $V_m$ have the form $V_m = (J_m, 
X_m, \tilde Z_m)$,
where $J_m \in \empJ$ is the state visited at time $T_m$ and $(\tilde Z_m, X_m)$ are covariates taking on value
in $E_1 = R^d \times [0,\tau], \tau< \infty$. 
The pair $(\tilde Z_m, X_m)$ may represent some
measurements taken upon entrance into the state $J_m$.
For any Borel set $B$ of $E_1$, let
$
\mu_{m+1}(B,t,j) = Pr ((\tilde Z_{m+1}, X_{m+1}) \in B| T_{m+1} = t, J_{m+1} = j, (T_l,J_l,\tilde Z_l,X_l)_{l=0}^m)
$
and suppose that
\beaa
& & Pr (T_{m+1}  -  T_m \leq s, J_{m+1} = j| (T_{\l}, J_{\l}, 
\tilde Z_{\l}, X_{\l})_{\l=0}^m) = \\
& & 1(J_m = i) \int_{[0, s]}  
\exp [ - \sum_{\l} \int_0^u e^{\beta^T Z_{i\l m}(v)}  \alpha_{i\l}
(v, X_{m}) dv] 
e^{\beta^T Z_{ijm}(u)} 
\alpha_{ij}(u,X_m) du \;, 
\eeaa
where $Z_{ijm}(u) = f_m (u, T_l, J_l, \tilde Z_l, X_l: l = 0, \ldots, m)$
is a fixed deterministic function $f_m$, left continuous in $u$. 
The process 
$
\tilde N_{ij}(t,B) = \sum_{m \geq 0} 1(T_{m+1}\leq t, J_{m+1} = j, J_{m} = i, 
(\tilde Z_{m+1},  X_{m+1}) \in B)
$
has  compensator given by
\beaa
\Lambda_{ij}(t,B) & = & \Lambda_{ij}(T_m, B) \\
& + & \int_{(T_m , t]}  
\mu_{m+1}(B, u, j) 1(J_m = i) e^{\beta^T Z_{ijm}(u-T_m)}\alpha_{ij}(u - T_{m}, X_m) du \;.
\eeaa
In particular, setting
 $B = E_1$ and using $\mu_{m+1}(E_1, T_{m+1}, j)1(T_{m+1} < \infty) = 1$, 
$$
\Lambda_{ij}(t) = \Lambda_{ij}(t, E) =
 \Lambda_{ij}(T_m) + \int_{(T_m , t]}  
1(J_m = i) e^{\beta^T Z_{ijm}(u-T_m)}\alpha_{ij}(u - T_{m}, X_m) du 
$$
is the  compensator of the counting  process 
$
\tilde N_{ij}(t) = \tilde N_{ij}(t,E) = 
\sum_{m \geq 0} 1(T_{m+1}\leq t, J_{m+1} = j, J_{m} = i) 
$,
registering  transitions among the adjacent states of the model.

In the following we assume
the random censorship model of Gill (1980). Thus the
times at which the process is observed is determined a process
$C(s) = \sum_{m \geq 1} 1(C_{m-1} < t \leq C_m)$,
where  
$0 = C_0 \leq  C_1  \leq \ldots \leq C_m...$ is an increasing sequence
such that 
$C_m \in [T_m,T_{m+1}]$
are stopping times with respect to the history $\{\F_t\}_{t \geq 0}$
and  $(C_m)_{m \geq 0}$ is conditionally
independent of $\{(T_m, J_m, \tilde Z_m, X_m)\}_{m \geq 0}$ given 
$(J_0, \tilde Z_0, X_0)$.
If $T_m = C_m$, then no information is
available on either the sojourn time $T_{m+1} - T_m$, the states
$(J_m, J_{m+1})$ or the covariates $(\tilde Z_m, X_m)$, \newline 
$(\tilde Z_{m+1}, X_{m+1})$. 
If $C_m = T_{m+1}$, then the sojourn time 
$T_{m+1}- T_m$, the adjoining
states $(J_m , J_{m+1})$ and the covariates  
$(\tilde Z_m, X_m), 
(\tilde Z_{m+1}, X_{m+1})$ 
are observable. Finally,  
if $T_m < C_m < T_{m+1}$, then the state $J_m$ and the covariates
$(\tilde Z_m, X_m)$ are visible while the
sojourn time $T_{m+1} - T_m$ is only known to exceed $C_m - T_m$.
We also assume that the censoring process is monotone
in the sense that 
$
T_m \leq C_m < T_{m+1} \implies C_{m'} = T_{m'} \quad \mbox{for all} \quad m' \geq m \;.
$
This condition stipulates that the process terminates once censoring
takes place.
To construct estimates of the unknown parameters, we use
a time transformation which replaces the chronological (or calendar) time
scale by the duration scale (Gill (1980), Oakes and Cui (1994)).  
For $m \geq 0$, let
\beaa
N_{ijm}(v) & = & 1( T_{m+1} - T_m \leq v, J_m = i, J_{m+1} = j, T_m= 
C_{m+1}) \;, \\
Y_{im}(v) & = & 1(T_{m+1} - T_m \geq v, C_m - T_m \geq v, J_m =i) \;,\\
M_{ijm}(v) & = & N_{ijm}(v) - \int_0^v Y_{im}(u) e^{\beta^T Z_{ijm}(u)}
\alpha(u, X_m)du \;.
\eeaa

\blem Suppose that $\{\phi_m(v), m \geq 0, v \geq 0\}$ is a sequence of left-continuous
random functions such that the process
$
\phi \circ L(t) = \sum_{m \geq 0} \phi_m(t - T_m) 1(T_m < t \leq T_{m+1})
$
is predictable with respect to the filtration $\{\F_t\}_{t \geq 0}$
and 
$
\empE \int_0^{\infty} [\phi \circ L]^2(s) \Lambda_{ij}(ds) < \infty.
$
Then
\beaa
\empE \sum_m \int_0^{\infty} \phi_m(u) N_{ijm}(du)  =  \empE 
\sum_m \int_0^{\infty} Y_m(u) \phi_m(u) 
e^{\beta^T_0 Z_{mij}(u)} \alpha_{ij}(u,X_m) du & & \\
\empE [\sum_m \int_0^{\infty} \phi_m(u) M_{ijm}(du)]^2  =  
\empE \sum_m \int_0^{\infty} Y_m(u) \phi^2_m(u) 
e^{\beta^T_0 Z_{mij}(u)} \alpha_{ij}(u,X_m) du & &
\eeaa
In addition, if $\{\phi_{1m}: m \geq 0\}$ and
$\{\phi_{2m}: m \geq 0\}$ are two such sequences, then
$$
\empE [\sum_m \int \phi_{1m}(u) M_{ijm}(du)][ \sum_m \int \phi_{2m}(u) M_{klm}(du)] 
= 0
$$
for pairs $(i,j) \not=(k,\l)$. 
\elem

Much in the same way as in Gill (1980),
this lemma follows from the Dominated Convergence Theorem, 
martingale properties
of the processes $\tilde M_{ij}$, and 
\beaa
\int_0^{\infty} [\phi \circ L](s) C(s)\tilde N_{ij}(s) & = &\sum_{m \geq 0} 
\int_0^{\infty}\phi_m(u) N_m(du) \\
\int_0^{\infty} [\phi \circ L]^k(s) C(s)\Lambda_{ij}(ds) & = & \sum_{m \geq 0} 
\int_0^{\infty}\phi_m(u)^k Y_m(u) e^{\beta^T_0 Z_{ijm}(u)} 
\alpha_{ij}(u, X_m) du \;.
\eeaa
The identies hold almost surely for $k = 1,2$. We omit the details.

\par

\setcounter{section}{3}
\setcounter{equation}{0}
\noindent{\bf 3. Estimation in single-type event  processes}

In this section  we assume that all events 
are of a single type. 
To estimate the baseline cumulative hazard function, we use
conditional Aalen-Nelson estimator (Beran (1981))
$$
\hat A(v;x, \beta) = {1 \over na}\int_0^v { N_i(du,x) \over
S^{(0)}_{-i}(u,\beta,x)} \;,
$$
where 
$
S^{(0)}_{-i}(u,\beta,x) = {1 \over (n-1)a} \sum_{j \not=i} 
S_j^{(0)}(u,\beta,x)
$
and for each $i = 1, \ldots, n$,
$$
N_i(u,x) = \sum_m K_n(x, X_{im}) N_{im}(u), \quad
S_i^{(0)}(u, \beta,x) = 
\sum_m Y_{im}(u) e^{\beta^T Z_{im}(u)} K_n(x, X_{im}) \;. 
$$
Here $K_n(x,w)$ is  the boundary kernel of 
M\"uller and Wang (1994),
\beaa
K_n(x,w)  & = &  
1(x-a \leq w \leq x+a) K_{11}({x - w \over a}) 
\quad \quad \quad \quad \tif a < x < \tau-a \\
 &= &   1( w \leq (1+q)a) K_{1q}(q - {w \over a}) 
\quad  \quad \quad \quad \quad \quad \tif 0 < x \leq a \\
& = &  1(\tau-a \leq w \leq (1+a)p) 
K_{p1}({\tau -w \over a} - p) 
\quad \tif \tau- a \leq x \leq \tau 
\eeaa
where
\beaa
& & K_{11}(r)  =  2 C(\mu) ({1 \over 2} )^{2\mu+2} 
 (1 + r)^{\mu} (1 - r)^{\mu}  \quad \quad \quad \; \;\mbox{(central region)} \\
& & K_{pq}(r)  =   \quad \quad \quad
  \quad 
\quad \quad \quad \quad \quad \quad \quad \quad \quad \; \; \;
\mbox{(left boundary region)} \\
& = &  C(\mu) ({1 \over p+q} )^{2\mu+2} 
 (p + r)^{\mu} (q - r)^{\mu-1}
[ 2 r ((p-q)\mu - q) + \mu(p-q)^2 + 2 q^2 ] \\
& & K_{pq}(r)  =  \quad \quad \quad
\quad  \quad \quad \quad \quad \quad \quad \quad \quad 
\quad \mbox{(right boundary region)} 
\\
&= &  C(\mu) ({1 \over p+q
 } )^{2\mu+2} 
 (p + r)^{\mu-1} (q - r)^{\mu}
[ 2 r ((p-q)\mu +p ) + \mu(p-q)^2 + 2 p^2 ] \;, 
\eeaa
 $p,q \in (0,1)$, and $C(\mu) = 2(2\mu+1) {2\mu - 1 \choose \mu}$.
The kernels $K_{pq}$ are Jacobi polynomials, and for 
$(p,q) = (1,1), (1,q)$ and $(p,1)$, we have
$$
\int_{-p}^q K_{pq}(u)du = 1, \int_{-p}^q u K_{pq}(u)du= 0,
\int_{-p}^q u^2 K_{pq}(u)du <\infty \;.
$$
Table 3.1 gives form of these kernels for polynomials of degree 2,4, and 6.

\begin{center}
Table 3.1 about here
\end{center} 
In the following we assume that  
$(u,x) \in \R = [0, \tau_0] \times [0, \tau], \tau_0 < \infty, \tau< \infty$. To control the bias of the risk process and the Aalen-Nelson estimator,
we  need the following regularity conditions.

\noindent
\bf Condition A \rm
\bitem
\item
[{(i)}] The variables $X_{im}$  have densities
$f_m(x)$  with respect to Lebesgue measure
on $[0, \tau]$.
\item[{(ii)}] There exists a bounded open
neighbourhood $\B$ of the true parameter
value $\beta_0$ such that 
 $\empE \sum_m 
[Z_{im}(u)]^{\otimes k}Y_{im}(u) \exp [\beta^T Z_{im}(u)] < \infty$, for
$k = 0,1,2$.
\item[{(iii)}] For $k = 0,1,2$, and $\beta \in \B$, the functions
$$
s^{(k)}(u,\beta,w) =  
\sum_m \empE (
[Z_{im}(u)]^{\otimes p}Y_{im}(u) \exp [\beta^T Z_{im}(u)]|X_m=w) f_m(w)
$$
are uniformly bounded and twice differentiable with respect to $\beta$.
In addition
$\nabla s^{(0)}(u,\beta,w) = s^{(1)}(u,\beta,w), \nabla^2 
s^{(0)}(u,\beta,w) = s^{(2)}(u,\beta,w)$,
and the functions $s^{(k)}(u,\beta,w), k = 0,1,2$ 
are uniformly Lipshitz continuous in
$\beta$.
\item[{(iv)}] The function $\alpha(u,w), (u,w) \in \R$ is bounded.
\item[{(v.1)}] The functions $s(u,\beta,w) = s^{(k)}(u,\beta,w), k = 0,1,2$,
and $\alpha(u,w)$ \newline 
satisfy $\sup\{|\alpha(u,w_1) - \alpha(u,w_2)|: (u,w_j),  \in \R,
|w_1 - w_2| \leq a,  j = 1,2\}  =  O(a)$ and
$\sup\{|s(u,\beta,w_1) - s(u,\beta, w_2)|: (u,w_j) \in \R,
|w_1 - w_2| \leq a, \beta \in \B, j = 1,2\} = O(a)$.

\item[{(v.2)}] $s(u,\beta,w)$ and $\alpha(u,w)$ are 
twice differentiable with respect to
$w$ with a uniformly bounded second derivatives $s''(u,\beta,w), 
\alpha''(u,w)$ such that \newline
$\sup\{|\alpha''(u,w_1) - \alpha''(u,w_2)|: (u,w_j),  \in \R,
|w_1 - w_2| \leq a,  j = 1,2\} =  O(1)$ and
$\sup\{|s''(u,\beta,w_1) - s''(u,\beta, w_2)|: (u,w_j) \in \R,
|w_1 - w_2| \leq a, \beta \in \B, j = 1,2\} = O(1)$.  
\eitem

We  refer to this condition as A.1 or A.2, depending on whether
the assumption (v.1) or (v.2) is in force. For $k = 1,2$, let
$S^{(k)}_{-i}(u, \beta,x) = \nabla^k S^{(0)}_{-i}(u, \beta,x)$
be the vector and matrix of first and second derivatives of the 
risk process $S_{-i}^{(0)}$ with respect to $\beta$.
Set     
 $\bar s^{(k)}(u,\beta,x) = a^{-1}\empE S^{(k)}_i(u,\beta,x)$,
$\bar n(u,x) = \empE N_i(u,x)$
and
$$
\bar A(v;x,\beta_0) = \int_0^v {\bar n(du;x) \over 
\bar s^{(0)}(u,\beta_0,x)} \;. 
$$

\bprop   Under assumptions A we have
$\bar s^{(k)}(u,\beta,x) - s^{(k)}(u,\beta,w) = O(a^r)$ for $k= 0,1,2$,
uniformly in
$(u, x) \in \R$ and $\beta \in \B$, 
and  $\bar A(v;x, \beta_0) - A_0(v;x) = O(a^r)$ uniformly in $(v,x) \in \R$. 
Here $r = 1$ under condition A.1 and 
$r = 2$ under condition A.2. 
\eprop

\Proof 
Dropping the superscript $k$, in the central region we  have
$$
 {1\over a}\empE S_i(u,\beta,x) = 
a^{-1} \int_{x-a}^{x+a} K_{11}({x-w \over a}) s(u,\beta,w)dw 
 =  \int_{-1}^1 K_{11}(r) s(u,\beta,x-ra)dr \;.
$$
In the left and right boundary regions, the expectation
$a^{-1} \empE S_i(u,\beta,x)$ is 
\beaa 
&  &  a^{-1}
\int_{x-qa}^{x+a}K_{1q}({x - w \over a}) s(u,\beta,w)dw
 =  \int_{-1}^q K_{1q}(r) s(u,\beta,x-ra) dr \\
& & a^{-1} \int_{x-a}^{x-pa} 
K_{p1}({x - w \over a}) s(u,\beta,w)dw 
 = \int_{-p}^1 K_{p1}(r) s(u,\beta,x-ra) dr  
\eeaa
In  the left boundary region, $q = x/a$ and in the right-boundary region
$p=(\tau-x)/a$.
Under condition (v.1), we have $|a^{-1}\empE S_i(u,\beta,x) - s(u,\beta,x)| = O(a)$,
uniformly in $(u,x) \in \R$ and $\beta \in \B$.
Under condition (v.2), we have
$$
a^{-1} \empE S_i(u,\beta,x) - s(u,\beta,x) = 
{a^2 \over 2}s''(u,\beta,x) \int_{-p}^q
r^2 K_{pq}(r) dr + O(a^2) \;.
$$
Similarly
$
\bar n(u,x)  = 
 \int_0^u \int_{-p}^q s^{(0)}(v, \beta_0, x-ra) \alpha(v, w-ra) 
K_{pq}(r) dr dv \;. 
$
Therefore, 
if  one of the two functions ($s$ or $\alpha$) is Lipschitz of order
1, then $\bar n(u,x) - \int_0^u s(v,\beta_0,x)\alpha(v,x)dv = 
O(a)$, whereas 
if both functions are twice differentiable in $x$, then the bias is 
$$
 {a^2\over 2}\int_0^u \bigg\{{\partial^2\over\partial x^2}
[s^{(0)}(v,\beta_0,x))\alpha(v,x)]\bigg\}dv \int_{-p}^q r^2 K_{pq}(r)dr +O(a^2)
$$
We also have
$
\bar A(v;x, \beta_0) - A(v;x) = 
\int_0^v \gamma(u,x) A(du;x) \;,
$
where $\gamma(u,x) = $ \newline $[\bar n(du,x)/
\bar s^{(0)}(u,\beta_0,x)\alpha(u,w)] - 1$. Thus the 
bias  is of order 
$O(a^r), r =1,2$ 
 \done

We turn now to estimation of the regression coefficients.
The first  method corresponds to an M-estimator
obtained by solving the score
equation $\tilde \Phi_n(\beta) = 0$, where
$$
\tilde \Phi_n(\beta)  =  {1 \over n} \sum_{i=1}^n \sum_m 
\int_0^{\tau_0} [Z_{im}(u) S^{(0)}_{-i}(u,\beta,
X_{im}) -  S^{(1)}_{-i}(u,\beta, X_{im})] N_{im}(du) \;.
$$
The analysis of this score equation requires 
only smoothness conditions A.1 and 
second moment bounds on the risk processes. 
For the sake of convenience, these moment
bounds are given in the appendix.  
Let
$$
V(u,\beta,x) = [{s^{(2)} \over s^{(0)}} - 
({s^{(1)} \over s^{(0)}} )^{\otimes 2}](u,\beta,x)
$$
\bprop
Suppose that the conditions A.1 and D.2 (i)--(ii) hold. 
Let
$\Sigma_1(\beta_0)  =  \int_{\R} ( V
 [s^{(0)}]^2 ) (u,\beta_0,x) \alpha(u,x)dudx$ and
$\Sigma_2(\beta_0)  =  $ \newline $\int_{\R} (V
 [s^{(0)}]^3 ) (u,\beta_0,x) \alpha(u,x)dudx$.
Suppose that  $\Sigma_1(\beta_0)$ is a non-singular matrix,
that  $na^2 \downarrow 0$ and $na \uparrow \infty$. 
With probability tending to 1, the score equation
$\tilde \Phi(\beta) = 0$ has a unique root $\tilde \beta$ 
and 
$\sqrt n (\tilde \beta - \beta_0)$ converges in distribution
to a mean zero normal variable with covariance 
$\Sigma^{-1}_1(\beta_0)
\Sigma_2(\beta_0))[\Sigma^{-1}_1(\beta_0)]^T$. 
\eprop

The proof is given in Appendix D.
The next Proposition deals with asymptotic normality of the 
Aalen-Nelson estimator. We 
need the following  consistency assumption on the risk function.

\noindent
\bf Condition B \rm
 Suppose that 
$\inf\{s^{(0)}(u, \beta,w): u \leq \tau_0, \beta \in \B, 
w \in [0 \vee x-a_n, x+a_n \wedge \tau]\} > 0$.
Moreover,  that under assumption A.$r, r = 1,2$, we have
$$
\max_i\empE\sup_{ \beta \in \B, u \leq \tau_0} |{S^{(0)}_{-i} - \bar s^{(0)}\over s^{(0)}}|
(u,\beta,x) \to 0
$$
for a bandwidth sequence $a = a_n \downarrow 0$ such that $na \uparrow \infty$
and $na^{2r+1} \downarrow 0$.

\bprop Suppose that  conditions A.$r (r = 1,2)$, B and D.1 
are satisfied.
For any root-n consistent
estimate 
$\hat \beta$ of the parameter $\beta_0$, the process 
$[\sqrt {na} [A(v;x,\hat \beta) - A(v;x)], v \leq \tau_0]$ converges weakly in 
$\l^{\infty}([0, \tau_0])$ to a mean zero Gaussian process $G(v,x)$ with 
covariance
$$
\cov [G(v,x), G(v',x)] = d_{p(x),q(x)}(K) \int_{[0, v \wedge v']} 
{A(du,x) \over 
s^{(0)}(u,\beta_0,x)} \;.
$$
Here $r =1$ under condition A.1 and $r = 2$ under assumptions of condition
A.2. Moreover,
 $d_{pq}(K) =\int_{-p}^q K^2_{pq}(w)dw$
and $p(x)=q(x) = 1$ if $a < x < \tau-a$, $p = 1, 
q(x) = a^{-1}x$ if $0 < x < a$ and
$p(x) = a^{-1}(\tau-x), q(x) = 1$ if $\tau -a < x < \tau$.
\eprop

Finally, we consider a partial score likelihood estimate of the regression
coefficient. It is obtained  by solving the 
the score equation
 $\Phi_n(\beta) = 0$, where
$$
\Phi_n(\beta) = 
{1 \over n} \sum_m \sum_{i=1}^n \int_0^{\tau_0} [Z_{im}(u) 
 -  {S^{(1)}_{-i} \over S^{(0)}_{-i}}(u,\beta, X_{im})] N_{im}(du) \;.
$$
Note that this score function is similar to that arising in 
the standard Cox regression, except that we use leave-one-out risk processes.
The choice of risk processes $S^{(k)} = \Sigma_{j=1}^n S_j^{(k)}, k = 1,2$,
is also possible. In both cases the resulting score functions 
form an approximate V process of degree 4 and  the difference between them
 converges in probability to 0, but only  
under stronger moment conditions than
those considered in the appendix D.

To analyze the score function $\Phi_n(\beta)$, we 
require condition A.2, moment conditions, and the following  
 uniform consistency assumption. 

\noindent
\bf Condition C \rm
Suppose that $\inf\{s^{(0)}(u, \beta,x):(u,x) \in \R, \beta \in \B\} >0$.
Moreover, that  
$$
\max_i \empE \sup_{ (u,x) \in \R, \beta \in \B} 
|{S^{(0)}_{-i} - \bar s^{(0)}\over s^{(0)}}|
(u,\beta,x)| \to 0
$$
for a bandwidth sequence $a_n \downarrow 0, na_n^2 \uparrow \infty, 
na^4_n \downarrow 0$.

\bprop
Suppose that conditions A.2, C, D.2 are
satisfied and  the matrix
$
\Sigma(\beta_0)  =  \int_{\R} (V
s^{(0)} )(u,\beta_0,x) \alpha(u,x)dudx 
$
 is  non-singular.
With probability tending to 1, the score equation
$\Phi_n(\beta) = 0$ has a unique root $\hat \beta$, and 
$\sqrt n (\hat\beta - \beta_0)$ converges in distribution
to a mean zero normal variable with covariance $\Sigma^{-1}(\beta_0)$.
\eprop

The proofs of these propositions are given in Appendices B-D.
Similar to the approach of Pons and Visser (2000) we  use
U-process theory. Whereas in their setting asymptotic normality results
for the estimate $\hat \beta$ were obtained based on analysis of U-statistics
of degree 2, in our case the term $R_{1n}$ of their  Proposition 3
satisfies
only $R_{1n} = \sqrt n O_p(1) 
\sup_{\beta, (u,x)} |S^{(0)} - \bar s^{(0)}|(u,\beta,x)$. 
(Here $S^{(0)}= (na)^{-1}\sum S_i^{(0)}$.) In the case 
of one jump processes with bounded time independent covariates, say,
results of Einmahl and Mason (2000) imply that the supremum is of order
$O(\sqrt {\log a^{-1}/na})$ a.s., 
so that the term $R_{1n}$ diverges to infinity. 
In the following we therefore use
expansions of higher order.

Except for moment bounds, the proofs of these propositions do not use
any special properties of the $Z$ process, and we do not require
uniform consistency of the derivatives  $S_{-i}^{(k)}, k = 1,2$.
On the other hand, assumptions
B and C require 
a more detailed specification of the covariate $Z$
in order to  apply
inequalities from empirical process theory. 
The following proposition gives one set of conditions under which 
these assumptions hold. We  consider the assumption C only.
Let $\R_{1n} =  \{(u,x) \in \R: a \leq x \leq \tau-a\}$,  
$\R_{2n} =  \{(u,x) \in \R: 0 < x \leq a\}$ and 
$\R_{3n} =  \{(u,x) \in \R: \tau-a < x \leq \tau\}$. 
Let $
\H_{pn} =  \{h(u,\beta,x): (u,x) \in \R_{pn}, \beta \in \B\}$, $p=1,2,3$,
where 
$
h(u,\beta,x) = \newline
s^{-1}(u,\beta,x)\sum_m Y_{m}(u) e^{\beta^T Z_{m}(u)} K_n(x,X_m)
$.
Note that for large n
$$
\max_i \empE \sup_{(u,x) \in \R_{pn} \atop 
\beta \in \B} |{S^{(0)}_{-i}-\bar s^{(0)} \over s^{(0)}}|(u,\beta,x)
$$
is of the same order as 
$\mu_{pn} = \empE \sup\{
|h-\empE h|(u,\beta,x): (u,x) \in \R_{pn},
\beta \in \B\}$.

\bprop
Suppose that for some $r > 2$ the bandwidth sequence satisfies
$a_n \downarrow 0$, $na_n \uparrow \infty$ $b_n = \log a_n^{-1}/ 
(na_n) \downarrow 0$, $a_n^{-1}b_n^{r/2-1} = O(1)$
and there exists a random variable $H_{1n}$, 
such that 
1) $\empE H_{1n}^r = O(1)$; \newline 
2) $ \Vert h(u,\beta,x) \Vert_{L_2(P)} 
\leq \sqrt a_n \Vert H_{1n} \Vert_{L_2(P)}$ and 3) 
$N_{[]}(\epsilon \Vert H_{1n} \Vert_{L_2(P)}, \H_{1n}, \Vert \cdot \Vert_
{L_2(P)}) \leq 
[A \epsilon^{-1}]^V$ for some finite constants $A$ and 
$V$ not depending on n and $\epsilon \in (0,1)$.
Then $\mu_{1n}=  O(\sqrt b_n)$.
If in addition there exist random variables $H_{pn}, p =2,3$, such that 
4) $\empE H_{pn}^2 = O(a)$ and 5)
$N_{[]}(\epsilon \Vert H_{pn} \Vert_{L_2(P)}, \H_{pn}, \Vert \cdot \Vert_
{L_2(P)}) \leq 
[A_p \epsilon^{-1}]^{V_p}$ for some finite constants $A_p$ and 
$V_p$ not depending on n and $\epsilon \in (0,1)$,
then in the boundary regions we have $\mu_{pn}= O((na_n)^{-1/2}), p = 2,3$.
\eprop

Here $\Vert \cdot \Vert_{L_2(P)}$ is the $L_2(P)$ norm, and 
$N_{[]}(\eta, \H_{pn}, \Vert \cdot \Vert_{L_2(P)})$ is the minimal number of brackets
of $L_2(P)$-size $\eta$ covering the class $\H_{pn}$.

\Proof
By Theorem 2.14.2 in van der Vaart and Wellner ((1996), p.240),
in the central region we have
\beq
\mu_{1n} 
 \leq  {1 \over a_n\sqrt n} J_{[]}(\sqrt a_n, \H_{1n}, \Vert \cdot \Vert_{L_2(P)}) \\
+ a^{-1}_n \empE H_{1n} 1(H_{1n} \geq \sqrt nc(\sqrt a_n))
\eeq
where
$
J_{[]} (\delta,\H,\Vert \cdot \Vert_{L_2(P)}) = 
\int_0^{\delta}  [1 + \log N_{[]}
(\epsilon \Vert H \Vert_{L_2(P)}, \H, \Vert \cdot \Vert_{L_2(P)})]^{1/2}
d\epsilon 
$
and
$c(\delta) = \delta \Vert H \Vert_{L_2(P)}/[1 + \log N_{[]}(
\delta \Vert H \Vert_{L_2(P)}, \H, \Vert \cdot \Vert_{L_2(P)})]^{-1/2}$.
For $\delta  = \sqrt {a_n}$
the first term of (3.1) is of order
$O(\sqrt b_n)$. Since 
$c(\sqrt {a_n}) =  O(\sqrt {a_n / \log a^{-1}_n})$, the second term is
bounded by
$a^{-1}_n
(\sqrt n c(\sqrt a_n))^{1-r} \empE H_{1n}^r 
 =  O(\sqrt b_n) 
O(a^{-1}_n b_n^{r/2-1}) = O(\sqrt b_n)$. The same theorem in
van der Vaart and Wellner (1996) implies that 
in the boundary regions we have
$\mu_{pn} 
= n^{-1/2}a^{-1}_n
O(J_{[]}(1, \H_{np},\Vert \cdot \Vert_{L_2(P)}) \Vert H \Vert_{L_2(P)}
= O((na_n)^{-1/2})$, $p = 2,3$ \done

Using a somewhat tedious argument, it is not difficult to show that 
conditions of this proposition are satisfied in the case of covariates
not dependent on u. Under added envelope conditions, the proposition
is also satisfied by Lipshitz continuous covariates, covariates
that form functions of bounded variation, etc.

\par

\setcounter{section}{4}
\setcounter{equation}{0}
\noindent
{\bf 4.  Multi-type event processes}

The results of the previous section extend  to the multistate setting
provided the state space of the process is ``small''. An example 
is provided by an illness-death process in which a person in ``healthy''
state (0) can either progress to a ``death'' state (2), or can first develop
a  reversible disease (state 1) and subsequently die. 
In the absence of censoring,
the  cumulative transitions
rates are given by
$$
\Lambda_{ij}(t) = \Lambda_{ij}(T_m) +1(J_m =i) 
\int_{(T_m,t]} e^{\beta^T Z_{ijm}(s - T_m)} \alpha_{ij}(s -T_m, X_m) ds
$$
for $t \in (T_m, T_{m+1}]$.
Similarly to multi-type processes in Andersen \it et al \rm (1993), estimation
of regression coefficients can be based on the score function
$$
\Phi_n(\beta) = {1 \over n} \sum_{i=1}^n \sum_h \sum_m
\int [Z_{ihm}(u) - {S_{-ih}^{(1)} \over S_{-ih}^{(0)}}(u, \beta, X_{im})]  
N_{ihm}(du) \;,
$$
where the sum extends over pairs $h = (0,1), (0,2), (1,2), (2,1)$ of possible
one-step transitions,
$$
S_{-ih}^{(0)}(u,\beta,x) = {1 \over a_h} \sum_{j\not=i} Y_{jhm}(u) 
e^{\beta^T Z_{jhm}(u)}K_n(x, X_{im}) \;,
$$
and  $S^{(1)}_{ih}$ is the derivative of this process with respect to
$\beta$. Note that the bandwidth sequence $a_h = a_{nh}$ is taken here 
to depend on the transition type $h$. The orthogonality relations
of Lemma 2.1 imply that the score function is asymptotically normal
with covariance matrix  $\sum_h \Sigma_h(\beta)$, where matrices 
$\Sigma_h$ assume a similar form as in Proposition 3.5.  
The M-estimator of Proposition 3.3 provides an alternative estimate.

Another example of a multi-type process is provided by progressive
multistate models. In this case  a subject may  move among a finite number
of transient states, but each such state can be visited at most once.
As an example of such a model
we consider data on 3020 bone marrow transplant
(BMT) recipients for acute myelogeneous leukemia (AML) and acute lymphoblastic
leukemia (ALL).
The
data were collected by the International Bone Marrow Transplant Registry
(IBMTR) during the period 1991-2000.  Only first
transplants in remission are  considered
and all patients received transplant  from an HLA-identical sibling.
Transplant recipients first receive  high doses 
of chemotherapy and radiation to destroy malignant cells in bone
marrow and
elsewhere. To rescue them from the toxicity of this 
therapy, they subsequently receive bone marrow 
cells from a suitably matched donor.  

In the following we donote by TX the transplant state.
It can be
followed by a number of complications, among them
graft--versus--host
disease (GVHD), relapse and death
in
remission. 
Two forms of GVHD are usually distinguished. Acute GVHD (AGVHD) occurs in
the first 2--3 months following transplant, whereas  chronic GVHD 
(CGVHD) occurs later
in time. 
We use time independent covariates corresponding to $X$= square root
of patient's age 
at transplant, and binary covariates represting  donor--recipient sex--match, 
($Z$), disease type 
and GVHD
 prophylaxis treatment. 
The square-root transformation of age
serves to reduce skewness of the data. 
Removal of T--cells from the
donor's bone 
marrow and posttransplant administration of immune supressive
drugs are the major
GVHD prophylactic treatments.

We are interested in the dependence of the intensities of one-step
transitions on age. 
In Figures 4.1-4.3 we show  plots of  the baseline
cumulative hazards $A_{ij}(v|x)$ as functions of $x$ . Note that
for fixed $x$, $A_{ij}(v|x)$ is an increasing function of $v$,
but for fixed $v$ this function may assume a variety of forms.
Figure 4.1  shows that cumulative hazards of transitions $\mbox{TX} \to 
\mbox{AGVHD}$, $\mbox{TX} \to 
\mbox{CGVHD}$ and $\mbox{AGVHD} \to 
\mbox{CGVHD}$ are increasing functions of age, and this monotonicity pattern
is most pronounced in the case of transitions into the CGVHD state.
The cumulative hazards of transitions  
$\mbox{TX} \to \mbox{death}$ and $\mbox{CGVHD} \to \mbox{death}$ are both
U-shaped functions, suggesting higher incidence of death among older and very 
young patients.  
Finally, the graphs of cumulative hazards of transitions into the 
relapse state are decreasing functions of age, though nearly constant in age
in the upper tail.  Note that in the case of transitions originating 
from the TX state, all 3020 subjects enter into the risk process. However,
transitions originating from the GVHD  states use only those
subjects who progress to the AGVHD and/or CGVHD state. In particular,
a total of 560 patients progressed into the CGVHD state. Subsequently
100 developed relapse and 170 died in remission. Thus  transitions
from the CGVHD  state are heavily censored. The relatively
small number of relapses accounts for the noisy graphs of the cumulative
hazards of the CGVHD $\to$ relapse state. 

\begin{center}
Figures 4.1--4.3 about here
\end{center}

The regression coefficients for the model are reported in Table 4.1.
As in any multistate analysis based on the proportional hazard
model, the regression coefficients do not have a clear
meaning. For example, male recipients receiving transplant from
a female donor are at higher risk for progression from the transplant
state into the AGVHD and CGVHD state, but are also at lower risk for 
direct (one-step) transition from the transplant into the relapse state.
The overall effect of this covariate on the occurrence of death in remission
or relapse cannot be, however, directly assessed based on regression
coefficients because patients who develop AGVHD are at higher risk
for death in remission, and also female-to-male transplant increases the
risk of CGVHD to relapse transition. Likewise, the direction of the 
regression coefficients corresponding to each of the GVHD prophylactic
treatments varies from one transition to another.  
Examples of parameters which can be used to summarize effects of covariates
on the occurrence of endpoint events were discussed in 
Klein, Keiding and Copelan (1993), Arjas and Eerola (1993) and Dabrowska, Sun and Horowitz
(1994). Their extension to the present setting is beyond the scope
of this paper.

\begin{center}
Table 4.1 about here
\end{center}

\par

\setcounter{section}{5}
\setcounter{equation}{0}
\noindent{\bf Appenidx A: Preliminaries}

Let $W_1, \ldots, W_n$ be iid random variables with some distribution P.
An (asymmetric) $U$ statistics of degree $m, m \geq 1$ 
is denoted by
$$
\empU_{n,m}(h) = {(n-m)! \over n!}
\sum_{(i_1, \ldots, i_m) \in I_n^m} h(W_{i_1}, \ldots, W_{i_m})
$$
where  $I_n^m$ is the collection of vectors $(i_1, \ldots,i_m)$ with distinct
coordinates, each in $\{1, \ldots,n\}$. Assuming that the kernel h satisfies
$\empE |h(W_{1}, \ldots, 
W_m) | < \infty$,  
the Hoeffding projection of degree m
of the kernel $h$ is denoted by 
$\pi_m h (W_1, \ldots, W_m)$. We have
$\pi_m h(W_1, \ldots, W_m) =  
\sum_{A \subset \{1, \ldots, m \} } (-1)^{m-|A|} 
\empE_A h(W_{1}, \ldots, W_m) \;,
$
where for $\emptyset \not= A = \{i_1, \ldots, i_p\}, 1 \leq p \leq m$, 
$\empE_A$ denotes conditional expectation
with respect to variables $\{W_{j}, j \in A\}$ and 
$\empE_{\emptyset} h(W_1, \ldots,W_m) = \empE h(W_1, \ldots,W_m)$.
 Then 
$\empU_{n,m}(\pi_m h)$ forms a canonical U statistics of degree m. 
For canonical U-processes indexed by classes of kernels changing with n,
Lemma 3.5.2, Remarks 3.5.4 and inequality (5.4.3) in de la Pe\~na and Gine
(1999) provide the following.

\blem 
Let $\{\empU_{n,m}(h): h \in \H_n\}$ 
be  a canonical U--process over a measurable class 
class $\H_n$ of (asymmetric) kernels of degree m.
If $\H_n$ forms a Euclidean class of functions for a square integrable
envelope
$H_n$, then 
$\empE n^{m/2}\Vert \empU_{n,m}(h) \Vert_{\H_n}$ \newline $= O( 
\empE [H_n(W_1, \ldots,W_m)^2]^{1/2})$.
\elem

A measurable class of functions $\H$ defined on some measure space
$(\Omega, \A)$ is Euclidean for envelope $H$ is $h \leq H$ for all
$h \in \H$, and there exist constants $A$ and $V$ such that 
$N(\epsilon \Vert H \Vert_{L_2(P)}, \H, \Vert \cdot \Vert_{L_2(P)}) \leq 
(A/\epsilon)^V$ for all $\epsilon \in (0,1)$ and all probability 
measures $P$ such that $\Vert H \Vert_{L_2(P)}  < \infty$
(Nolan and Pollard, 1987). 
Here $\Vert \cdot \Vert_{L_2(P)}$ is the $L_2(P)$ norm and 
$N(\eta, \H, \Vert \cdot \Vert_{L_2(P)})$ is the minimal number
of $L_2(P)$--bals of radius $\eta$ covering the class $\H$.
In the case of classes $\H_n$ changing with $n$,
the Euclidean constants $A$ and $V$ are taken to be independent of n.

In the following we shall use U processes of degree $m \leq 1, 2, 3,4$.
Finally, in our case for each subject $i$, 
the sequence $W_i$ represents the total number of events
observed in the interval $[0, \tau_0]$, their times of the occurrence,
types and covariates observed at each jump time. 
The Euclidean property of the classes of functions appearing in the 
remainder of the text can be easily verified based on 
results of 
Nolan and Pollard (1987), 
Pakes and Pollard (1989) and Gin\'e and Guillou (1999).

\par

\setcounter{section}{6}
\setcounter{equation}{0}
\noindent{\bf Appendix B: Regularity conditions and two lemmas}

We give  some additional regularity conditions.

\noindent
\bf Conidtion D.0 \rm
(i) For sequences $(m) =(m_1, m_2), m_1 \not=m_2$ 
of nonnegative integers the variables $X_{(m)} = (X_{m_1}, X_{m_2})$ 
have joint density $f_{(m)}$ with respect to Lebesgue measure
on $[0, \tau]^2$.

(ii)
For sequences 
 $[m] =(m_1,m_2,m_3)$ of distinct
nonnegative integers, 
the variables   $X_{[m]} =
(X_{m_1}, X_{m_2}, X_{m_3})$ have joint densities
$f_{[m]}$ with respect to Lebesgue measure on 
 $[0, \tau]^3$.

For any vector,  we denote  by $| \cdot |$ the $\l_1$ norm. Without
loss of generality we assume that the neighbourhood $\B$ surrounding the 
true parameter $\beta_0$ corresponds to a ball $\B = \{\beta: |\beta -
\beta_0| \leq c_B\}$.

For nonnegative integers $p$ and $m$ define
$
\bar \theta^{(p)}_m(u) = |Z_m(u)|^p Y_m(u) e^{[|\beta_0| + c_B]|Z_m(u)|}$
and
$\theta^{(p)}_m(u,\beta)  = |Z_m(u)|^p Y_m(u) e^{\beta^T Z_m(u)}$.
For
$u \in [0, \tau_0]$, $\underline u = (u_1, u_2) \in [0, \tau_0]^2,
\bar u = (u_1, u_2, u_3) \in [0, \tau_0]^3$, and
$w \in [0, \tau]$, $\underline w = (w_1, w_2) \in [0, \tau]^2,
\bar w = (w_1, w_2, w_3) \in [0, \tau]^3$, let
\beaa
\sigma_{p_1, p_2}(\underline u, w) & = & \sum_m 
\empE [\prod_{j=1}^2 \bar \theta_m^{(p_j)} 
(u_j)|X_m = w] f_m(w) \;,\\
\rho_{p_1, p_2}(\underline u, \underline w) &
= &
\sum_{(m)} \empE [\prod_{j=1}^2 \bar \theta_m^{(p_j)}
(u_j)|X_{(m)} = \underline w] f_{(m)}(\underline w) \;,\\
\kappa_{1;p}(\bar u,w) & = &
\sum_{ m } \empE[
\theta_m^{(p)}(u_1, \beta_0) \prod_{j=2}^3 \theta_m^{(0)}
(u_j, \beta_0)|X_m=w] f_m(w) \;,\\
\kappa_{2;p}(\bar u,\underline   
w) & = &
\sum_{(m)} \empE[
\theta_{m_1}^{(p)}(u_3, \beta_0) \prod_{j=1}^2 \theta_{m_j}^{(0)}
(u_j, \beta_0)|X_{(m)} = \underline
w] f_{(m)}(\underline w) \;,\\
\kappa_{3;p}(\bar u,\bar
w) & = &
\sum_{[m]} \empE[
\theta_{m_1}^{(p)}(u_1, \beta_0) \prod_{j=2}^3 \theta_{m_j}^{(0)}
(u_j, \beta_0)|X_{[m]} = \bar
w] f_{[m]}(\bar w) \;,\\
s_{0;2}(u,\underline w) & = & 
\sum_{(m)} \empE[
                    \theta_{m_1}^{(0)}
(u, \beta_0  )|X_{(m)} = \underline
w] f_{(m)}(\underline w) \;, \\
s_{0;3}(u,\bar       w) & = & 
\sum_{[m]} \empE[
                   \theta_{m_1}^{(0)}
(u, \beta_0  )|X_{[m]} = \bar
w] f_{[m]}(\bar       w) \;.
\eeaa
Under conditions D.1 and D.2 these expectations exist, at least
in local neighbourhoods of a point $x \in [0, \tau]$.
Such local neighbourhoods correspond to sets
$
\R(x) = \{(u,w) \in \R: |w - x| \leq a\}
$.

\noindent
\bf Conidtion D.1 \rm
(i) The condition D.0 (i) is satisfied and
for integers
$p_1, p_2$ such that $p_j \geq 0, p_1+p_2 \leq 4$, we have
\beaa
& & \sup\{\sigma_{p_1, p_2}(\underline u,w):(u_1,w) \in \R(x)
, (u_2,w) \in \R(x) \} = O(1) \;.\\
& & \sup\{|\rho_{p_1,p_2}(\underline u, \underline w)
  )|: (u_j,w_j) \in \R(x),
                    j = 1,2 \} = O(1) \;.
\eeaa
(ii) The condition D.0 (ii) is sastisfied, and
\beaa
& & \sup \{ \kappa_{1;0}(\bar u,w):(u_j, w) \in \R(x), j = 1,2,3\} = O(1) \;,\\
& & \sup \{ \kappa_{2;0}(\bar u,\underline w):(u_1, w_1) \in \R(x), 
(u_2,w_2) \in \R(x), (u_3,w_1) \in \R(x)\} = O(1) \;,\\
& & \sup \{ \kappa_{3;0}(\bar u,\bar w):(u_j, w_j) \in \R(x), j = 1,2,3\} 
= O(1) \;,\\
& & \sup \{ s_{0;2}( u,\underline w):(u, w_j) \in \R(x), j = 1,2\} 
= O(1) \;,\\
& & \sup \{ s_{0;3}(u,\bar w):(u, w_j) \in \R(x), j = 1,2,3\} = O(1) \;.
\eeaa

\noindent
\bf Condition D.2 \rm
(i) The condition D.0 (i) is satisfied
and, for integers
$p_1, p_2$ such that $p_j \geq 0, p_1+p_2 \leq 4$, we have
\beaa
& & \sup\{\sigma_{p_1, p_2}(\underline u,w):(u_1,w) \in \R, (u_2,w) \in \R\}
 = O(1) \;,\\
& & \sup\{|\rho_{p_1,p_2}(\underline u, \underline w)
  )|: (u_j,w_j) \in \R,
|w_2 - w_1| \leq a, j = 1,2 \} = O(1) \;.
\eeaa
(ii) The condition D.0 (ii) is satisfied and,
for
$p=0,1$,  we have
\beaa
& & \sup \{ \kappa_{1;p}(\bar u,w):(u_j, w) \in \R, j = 1,2,3\} = O(1) \;,\\
& & \sup \{ \kappa_{2;p}(\bar u,\underline w):(u_1, w_1) \in \R, 
(u_2,w_2) \in \R, (u_3,w_1) \in \R, |w_2 - w_1| \leq a\} = O(1) \;,\\
& & \sup \{ \kappa_{3;p}(\bar u,\bar w):(u_j, w_j) \in \R,
|w_2 - w_1| \leq a, |w_3 - w_2|, j = 1,2,3\} = O(1) \;.
\eeaa

We now give two lemmas which collect bounds on certain
random variables arising in the analysis of the Aalen-Nelson
estimate. Both  can be verfied using elementary algebra,
H\"older's inequality and conditions A and D.

\blem
Suppose that 
$\inf \{s^{(0)}(u,\beta,x): \beta \in \B, u \leq \tau_0\} > 0$. 
For $k = 0,1,2$, let
$
\bar f_{kni}(u, \beta,x)  =  [s^{(0)}(u, \beta,x)]^{-1} 
\sum_m \theta^{(k)}(u, \beta) |K_n(x, X_{im})|$ and
$f_{kni}^*(u, \beta,x)  =  [s^{(0)}(u, \beta,x)]^{-1}
\sum_m \theta^{(k)}(u, \beta) \alpha(u,X_{im})|K_n(x, X_{im})|$.
If conditions A and D.1 (i) hold, 
then
$a^{-1} \empE \prod_{p=1}^2 \bar f_{k_pni}(u_p, \beta,x) = O(1)$ and
$a^{-1} \empE \prod_{p=1}^2 f_{k_pni}^*(u_p, \beta,x)
 = O(1)$, uniformly in $u_1, u_2 \leq \tau_0$ and $\beta \in \B$.
If in addition the condition D.1 (ii) holds, then
$a^{-1} \empE \prod_{p=1}^3 f_{k_pni}^*(u_p, \beta_0,x) = O(1)
$ uniformly in $u_1, u_2, u_3 \leq \tau_0$.
If $\inf \{s^{(0)}(u,\beta,x): \beta \in \B, (x,u) \in \R\} > 0$
and conditions D.2 hold, 
then these bounds are also uniform in $x, x \in [0, \tau]$.
\elem

\blem
Supose that $\inf \{s^{(0)}(u,\beta_0,w): u \leq \tau, \beta \in \B, 
w \in [x -a_n \vee0, x+a_n \wedge \tau]\} > 0$. 
Set
\beaa
{\bar {\bar f}}_{ni}(u,x)  & = & [\bar s^{(0)}(u,x)]^{-2} 
[\bar S^{(2)}_i(u,x) + 
\bar S^{(1)}_i(u,x) \bar s^{(1)}(u,x)] \;.\\
\bar S^{(p)}_i(u,x)  & = &   
\sum_m \bar \theta^{(p)}_{im}(u) 
|K_n(x,X_{im})| \;,\\
\bar s^{(0)}(u,x)   & = &  \sum_m \empE 
Y_{im}(u)\exp ([-|\beta_0|-c_B] |Z_{im}(u))| X_{im} =x) f_m(x) \;,\\
 \bar s^{(1)}(u,x)  & = &       \sum_m \empE (\bar \theta^{(p)}_m
| X_{im} =x) f_m(x) \;,\\
{\bar {\bar g}}_{nj}(v,x) & = & 
\sum_m \int_0^v |K_n(x, X_{jm})| [\bar s^{(0)}(u,x)]^{-1}
N_{jm}(du) \;,\\
\bar g_{in}(u, \beta,x) & = & \sum_m
\int_0^u |K_n(x,X_{im})| [ s^{(0)}(u,\beta,x)]^{-1}
N_{im}(du) \;,\\ 
\eeaa
and let
\beaa
H_{0n}(W_i) & = & a^{-1/2}  [\bar g_{ni}(\tau_0, \beta_0,x)
+ \int_0^{\tau_0} f_{0ni}^*(u, \beta_0,x) du] \;,\\
H_{1n}(W_i) &=& a^{-1/2} \sum_{m} \int_0^{\tau_0} {|K_n(x,X_{im})| \over s^{(0)}(u,\beta_0,x)} Y_{im}(u)
e^{\beta^T_0 Z_{im}(u)} \quad \times \\
& & \quad \quad \quad \times \quad |\alpha(u, X_{im}) - \alpha(u,x)| du \;,\\
H_{2n}(W_i,W_j) & = & { 1 \over a \sqrt {na}} \int_0^{\tau_0} 
\bar f_{0ni}(u, \beta_0,x) \bar g_{nj}(u,\beta_0,x) \;,\\
H_{3n}(W_i) & = &  a^{-1/2} 
\int_0^{\tau_0} {|\bar s^{(0)} - s^{(0)}| \over s^{(0)}} (u,\beta_0,x) 
\bar g_{ni}(du, \beta_0,x) \;,\\
H_{4n}(W_i) & = & a^{-1/2}  \int_0^{\tau_0}
\bar f_{ni}(u, \beta_0,x) 
 | \alpha(u,x)du -  {\empE N(du,x) \over s^{(0)}(u,\beta_0,x)}| \;,\\
H_{5n}(W_i,W_j) & = & {1 \over na^2}
\int_0^{\tau_0} \bar f_{1ni}(u,\beta_0,x) \bar g_{nj}(du,\beta_0,x) 
\;,\\
& + &  {c_B \over na^2}
\int_0^{\tau_0} {\bar {\bar f}}_{ni}(u,x) {\bar {\bar g}}_{nj}(du,x) \;.\\
\eeaa
If conditions  A.$r (r=1,2)$ and  D.1  hold, then
$\empE H_{0n}^2(W_1)  =  O(1)$, $\empE H_{0n}^3(W_1) = O(a^{-1/2})$, 
$\empE H_{1n}^2(W_1)  =  O(a^2)$, 
 $\empE H_{3n}(W_1)^2  =   O(a^{2r})$ and  
$\empE H_{4n}(W_1)^2 =   O(a^{2r})$  We also have 
$\empE H_{2n}^2(W_1,W_2)   =   O(({na})^{-1})$, 
$\empE H_{5n}^2(W_1,W_2)   =   O(({na})^{-2})$ and
$n\empE [\empE_{\{1\}} H_{5n}(W_1,W_2)]^2   =   O(({na})^{-1})
= n\empE [\empE_{\{2\}} H_{5n}(W_1,W_2)]^2$.
\elem

\par

\setcounter{section}{6}
\setcounter{equation}{0}
\noindent{\bf Appendix C: 
Proof of Proposition 3.4} \rm

Set
$$
b(v,x) = \int_0^v {\alpha(u,x) \over s^{(0)}(u,\beta_0,x)} 
 [{\gamma(u,x) - \bar s^{(0)}(u,\beta_0,x)} ]du \;,
$$
where $\bar s^{(0)} (u,\beta_0,x) = \empE S^{(0)}_{-i}(u,\beta_0,x)$, 
$\gamma(u,x) = { \bar n(du,x) / \alpha(u,x)}$
and $\bar n(v,x) = \newline \Sigma_m \empE  N_{im}(v)K_n(x,X_{im})$. Then 
$\sqrt {na} [\hat A(v;x, \beta_0) - A_0(v,x) - b(v,x) ]= \hat Z_n(v,x) +  
R_n(v,x)$,
where 
$$
\hat Z_n(v,x) = 
 \sqrt {n \over  a} \sum_{i=1}^n \sum_{m} \int_0^v 
{K_n(x, X_{im}) \over s^{(0)}(u,\beta_0,x)} M_{im}(du) + R_n(v,x) \;,
$$
and $R_n(v,x)$ is a remainder term given below.
Under conditions A.$r, r = 1,2$, we have 
$\sqrt {na} b(v,x) = O(\sqrt {na} a^r) = o(1)$. Therefore
it is enough to  show  
 that the process 
$\hat Z_n(v,x)$  converges in $\l^{\infty}([0, \tau_0])$
to a time transformed Brownian motion and the remainder term $R_n$
is asymptotically negligible.
 
We have
$\hat Z_n(v,x) = \sqrt n [P_n-P] h_{n,v}$ where
$$
h_{n,v}(W_i) = a^{-1/2} \sum_{m} \int_0^v 
{K_n(x, X_{im}) \over s^{(0)}(u,\beta_0,x)} M_{im}(du) \;.
$$ 
The class $\H_n = \{h_{n,v}: v \leq \tau_0\}$ 
consists of functions that can be represented
as a linear combination of at most four monotone functions with respect to
$v$ and has   envelope $4 H_{0n}(W_i)$. 
By Lemma 6.9 we have (i) $\empE H_{0n}^2(W_1) = O(1)$ 
and (ii) $
\empE H_{0n}(W_1) 1(H_{0n}(W_1) > \eta \sqrt n)  \leq   
\empE H_{0n}^3(W_1)( \eta 
\sqrt {na})^{-1} \to 0$ for any $\eta > 0$.
Also (iii) for any $0 < v_1 < v_2 \leq \tau_0$, 
the difference $
|h_{nv_1} - h_{nv_2}|(W_i)$ is bounded by
$$ a^{-1/2} \sum_{m}
\int_{v_1}^{v_2} {|K_n(x,X_{im})| \over s^{(0)}(u,\beta_0,x)} M_{im}(du) 
 +   {2 \over \sqrt a}
\int_{v_1}^{v_2} f^*_{0ni}(u, \beta_0,x) du \;.
$$
Using $(x+y)^2 \leq 2(x^2 + y^2)$ and Lemma 2.1, 
$\empE |h_{nv_1} - h_{nv_2}|^2 (W_1)$ is bounded by
$$  
{2 \over a} \int_{v_1}^{v_2} \int_0^{\tau}
{[s^{(0)} \alpha](u,\beta_0,w)
\over s^{(0)}(u,\beta_0,x)^2} 
K^2_n(x,w) dw du + {8 \over a} \int_{v_1}^{v_2} \int_{v_1}^{v_2} 
\empE \prod_{p=1}^2 f_{0ni}^*(u_p, \beta_0,x) du_1 du_2 \;, 
$$
and is  of order  $O(|v_2 - v_1| + |v_2 - v_1]^2)$. 
Lemmas 2.1 and 3.2, imply (iv) 
\beaa
\var [\hat Z_n(v_1,x)]  = d_{p(x),q(x)}(K) 
\int_0^{v_1} {\alpha(u,x) \over s^{(0)}(u, x, \beta_0)} du
+ O(a) \;,
\eeaa
and $cov [\hat Z_n(v_1,x), \hat Z_n(v_2,x) -  \hat Z_n(v_1,x)] = O(a)$.
Finally,  (v) the class of functions 
$\{h_{nv}: v \leq\tau_0\}$ has polynomial bracketing
number. Properties (i)-(v) and 
Theorem 2.11.23 in van der Vaart 
and Wellner (1996) imply that 
$\{\hat Z_n(v,x): v \in \tau_0\}$ converges weakly 
$\l^{\infty}([0, \tau_0])$  to a tight
Gaussian process.

The remainder term $R_n(v,x)$ is given by
$R_n(v,x) = \sum_{j=1}^5 R_{jn}(v,x)$ where
\beaa
R_{1n}(v,x) & = & { 1 \over \sqrt {na}}  \sum_{i=1}^n 
\int_0^v \tilde f_{ni}(u,x) du - \sqrt {na} b(v,x) \;,\\
R_{2n}(v,x) & = & -{\sqrt {na} \over n(n-1) a^2} \sum_{i \not=j} 
\int_0^v [f_{ni} - \empE f_{ni}](u,x) [
g_{nj} - \empE g_{nj}](du,x)] \;,\\
R_{3n}(v,x) & = & - {1 \over \sqrt {na}} \sum_{i=1}^n  \int_0^v 
[{\bar s^{(0)} - s^{(0)} \over s^{(0)}}] 
(u,\beta_0,x) [g_{ni}-\empE g_{ni}](du,x) \;,\\
R_{4n}(v,x) & = & {1 \over \sqrt {na}} \sum_{i=1}^n\int_0^v ([f_{ni} -
\empE f_{ni}] s^{(0)}) (u,\beta_0,x) \alpha(u,x)du - 
\empE N(du,x)] \;,\\
R_{5n}(v,x) & = & {1 \over \sqrt {na}}  \int_0^v 
[{\bar s^{(0)} -s^{(0)} }](u,\beta_0,x) \left [ 
\alpha(u,x)du - {\empE N(du,x) \over s^{(0)}(u, \beta_0,x)} \right] \;,\\
R_{6n}(v,x) & = & \sqrt {na} \sum_{i=1}^n \int_0^v 
{[S^{(0)}_{-i}(u,\beta_0,x) - s^{(0)}(u,\beta_0,x)]^2 \over 
S^{(0)}_{-i}(u,\beta_0,x) s^2(u,\beta_0,x)} \bar N_i(du,x) \;,
\eeaa
where
\beaa
f_{ni}(u,x) & = & [ s^{(0)}(u, \beta_0,x)]^{-1}
\sum_m Y_{im}(u) e^{\beta_0^T Z_{im}(u)} K_n(x,X_{im}) \;,\\
\tilde  f_{ni}(u,x) & =  &  [s^{(0)}(u, \beta_0,x)]^{-1} \sum_m 
[\alpha(u,X_{im}) - \alpha(u,x)]Y_{im}(u) e^{\beta_0^T 
Z_{im}(u)} K_n(x,X_{im}) \;,\\
g_{ni}(u,x)  & = &  
\sum_m\int_0^u K_n(x,X_{im}) [s^{(0)}(v,\beta_0,x)]^{-1} N_{im}(dv) \;.
\eeaa

The term $R_{1n}$
has mean zero. By decomposing the 
integrands and the integrators into their positive
and negative parts, we have $(na)^{-1/2}R_{1n}(v,x) + b(v,x) = \empP_n h_{1nv}$ where $h_{1nv}(W_i)$ is a sum
of four monotone functions, bounded by $H_{1n}(W_i)$.
Thus $R_{1n}(v,x)$ is a normalized empirical process over a Euclidean
class of functions for envelope $ 4 H_{1n}(W_i)$.
By Lemmas 6.9 and 5.7, we have $\empE H_{1n}(W_1)^2 = O(a^2)$
and $\empE\sup_v |R_{1n}(v,x)| =
O(a)$. Similarly, using envelopes $H_{3n}$ and $H_{4n}$, 
we can show that 
$ \empE\sup_v |R_{3n}(v,x)| = O(a^r) = 
\empE\sup_v |R_{4n}(v,x)|$ 
and
$R_{5n}(v,x) = O(\sqrt {na} a^{2r})$ a.s.uniformly in $v \leq\tau_0$.
The term 
 term $R_{2n}$ is easily seen to form  a canonical U-process of degree 2
over a Euclidean class of functions
with  envelope  $H_{2n}'(W_i,W_j) =H_{2n}(W_i,W_j) + \empE_{\{1\}}
H_{2n}(W_i,W_j) + \empE_{\{2\}} 
H_{2n}(W_i,W_j) +  \empE H_{2n}(W_i,W_j)$. Lemmas 6.9 and 5.7
imply
$ \empE\sup_v |R_{2n}(v,x)|$ \newline $=
O((na)^{-1/2}))$, since
$\empE H^2_{2n}(W_i,W_j) =
O((na)^{-1})$ and  $\empE [H_{2n}']^2(W_1,W_2)$ is of the same order.

Next define
\beaa
R_{7n} & = & {1 \over \sqrt {na}} \sum_{i=1}^n
\int_0^{\tau_0} ({S_{-i}^{(0)} - s^{(0)} \over 
s^{(0)}} )^2 (u, \beta_0,x)
 \bar g_{ni}(du,\beta_0,x) \\
& \leq &  2 \sqrt {na} O(a^{2r}) {1 \over na} \sum_{i=1}^n \
\bar g_{ni}(\tau_0, 
\beta_0,x) + O(1) (R_{7n;1} + R_{7n;2}) \;,
\eeaa
where 
$R_{7n;1} =
\sqrt {na}  a^{-3} \empU_{n,3}(h)$, $R_{7n;2} =\sqrt {na} (na^3)^{-1}
\empU_{n,2}(\bar h)$,
$$
h(W_i,W_j,W_k)  =  \int_0^{\tau_0} [(f_{0nj} - \empE f_{0nj}) 
(f_{0nk} - \empE f_{0nk}) ](u,\beta_0,x) \bar g_{ni}(du,\beta_0,x) 
$$
and $\bar h(W_i,W_j) = h(W_i,W_j,W_j)$.
The first term is of order $O_p(\sqrt {na} a^{2r})$. 
We have $\empE H(W_1,W_2,W_3) = 0$ and, 
using Lemmas 6.9 and   5.7, 
$\empE (na)^{1/2}a^{-3} |
\empU_{n,3}(\pi_3 h)|  = 
O((na)^{-1})$ and
$\empE(na)^{1/2}a^{-2} | \empU_{n,2}(
\pi_2 [\empE_{\{23\}} h])| = O((na)^{-1/2})$. The remaining 
projections are 0.
In the case of the term $R_{7n;2}$, we have $\empE R_{7n;2} = O((na)^{-1/2})$
and the expected $\empE |R_{7n;2}|$ is of the same order.

We  consider now  term $R_{6n}$. 
For $\epsilon \in (0,1)$, define
$$
\Omega_n(\epsilon) =  \bigl \{ {1 \over 1+\epsilon} \leq \min_i \inf_{u \leq \tau_0 
\atop \beta \in \B}
{s^{(0)} \over S^{(0)}_{-i}}(u,\beta,x) \leq  
\max_i \sup_{u \leq \tau_0 \atop
\beta \in \B} {s^{(0)} \over S^{(0)}_{-i}}(u,\beta,x) \leq  
{1 \over 1 - \epsilon} \bigr \} \;.
$$
We have
$P(\Omega_n(\epsilon)) \leq \min_i P(
\sup_{u \leq \tau}
|S^{(0)}_{-i}/s^{(0)} -1 |(u, \beta,x)| \leq \epsilon) \to 1$, by condition
B and Markov's inequality. 
On the event $\Omega_n(\epsilon)$,
we also have \newline
$\sup_{v \leq\tau_0} |R_{6n}(v,x)| \leq (1- \epsilon)^{-1} R_{7n}$.   
Therefore 
$P(\sup_{v \leq\tau_0}|R_{6n}(v,x)| > \eta) \leq P(\Omega_n^c(\epsilon)) 
+ 
P(\sup_{v \leq\tau_0}|R_{6n}(v,x)| > \eta, \Omega_n(\epsilon)) \leq
P(\Omega_n^c(\epsilon)) +
P(R_{7n} > (1-\epsilon)\eta) \to 0$ for any $\eta > 0$.

Finally, suppose that $\hat \beta$ is a $\sqrt n$ consistent estimate of the 
parameter $\beta$. Then
\bea
& & \sqrt {na}[\hat A(v,x, \hat \beta) - \hat A(v,x, \beta_0)] = \nonumber \\
& &  \sqrt n[\hat \beta - \beta_0] \sqrt a \int_0^v {S^{(1)}_{-i} \over 
[S^{(0)}_{-i}]^2}(u, \beta^*,x) N_i(du,x) \;,  
\eea
where $\beta^*$ is  between $\beta_0$ and $\hat \beta$.
Let $I_n(\beta) = a^{-2} U_{n2}(h_{\beta})$, where
$h_{\beta}(W_i,W_j) = \int_0^{\tau_0}
\bar f_{1ni}(u,\beta,x) \bar
g_{nj}(du,\beta,x)$.
It is easy to see that $\empE I_n(\beta) = O(1)$. 
By Lipschitz continuity of the function $h_{\beta}$
with respect to $\beta$,
$I_n(\beta)$ is  a U-process of degree 2 over a Eulidean class of functions
for envelope $H_{5n}(W_i,W_j)$. 
By Lemmas 6.9 and  5.7,
$\empE \sup_{\beta \in \B}| a^{-2} \empU_{n,2}(\pi_2 h_{\beta})| =
O ([\empE H_{5n}(W_1,W_2)^2]^{1/2}) =
O((na)^{-1})$,
$\empE \sup_{\beta \in \B} |a^{-2} U_{1n} (\pi_1 \empE_{\{1\}} h_{\beta})| 
= O((na)^{-1/2}))$, $
\empE \sup_{\beta \in \B} |a^{-2} U_{1n} (\pi_1 \empE_{\{2\}} h_{\beta})|
$ $ =  O((na)^{-1/2}))$.
Therefore $\sup_{\beta \in \B} |I_n(\beta)| = O_p(1)$. 
Further, if $\hat \beta$ is a $\sqrt n$ consistent estimate of $\beta_0$,
then $\sqrt n [\hat \beta - \beta_0] = O_p(1)$. 
To show that the right--hand side of (6.1) is of order $O_p(\sqrt a)$,
it is enough to note that for any $\epsilon \in (0,1)$,   
the supremum $\sup \{ |\int_0^v
S^{(1)}_{-i}  [S_{-i}^{(0)}]^{-2}(u,\beta,x) N_i(du,x) |: v \leq\tau_0,
\beta \in \B \}$ 
is bounded by 
$(1- \epsilon)^{-2} \sup_{\beta \in \B} 
I_n(\beta)$ on the event $\Omega(\epsilon)$.

\par

\setcounter{section}{7}
\setcounter{equation}{0}
\noindent{\bf Appendix D: Proof of Propositions 3.3 and 3.5}

Define
\beaa
\tilde \Phi_{0n}(\beta_0) & = & {1 \over \sqrt n} \sum_{i=1}^n \sum_m \int
[Z_{im}(u) s_{-i}^{(0)}(u, \beta_0, X_{im}) - 
s_{-i}^{(1)}(u, \beta_0, X_{im})] M_{im}(du) \;,\\ 
\tilde \Sigma_{n}(\beta) & = & {1 \over  n} \sum_{i=1}^n \sum_m \int
[S^{(2)}(u, \beta, X_{im}) -
Z_{im}(u) \otimes S^{(1)}(u, \beta, X_{im})] N_{im}(du) \;,\\
\Phi_{0n}(\beta_0) & = & 
{1 \over \sqrt n} \sum_{i=1}^n \sum_m \int
[Z_{im}(u) -
{s^{(1)} \over s^{(0)}}
(u, \beta_0, X_{im})] M_{im}(du) \;,\\
\Sigma_{0n}(\beta) & = & {1 \over n} \sum_{i=1}^n \sum_m \int [{s^{(2)} \over 
s^{(0)}} - \left (s^{(1)} \over s^{(0)} \right )^{\otimes 2}](u, \beta, 
X_{im}) N_{im}(du) \;,\\
\Phi_{1n}(\beta_0) & = & {1 \over  \sqrt n}  \sum_m
\int_0^{\tau_0} 
[Z_{im}(u){S_{-i}^{(0)} \over s^{(0)}}(u, \beta_0, X_{im}) 
 - {S^{(1)}_{-i} \over s^{(0)}}(u, \beta_0, X_{im}) 
]
 N_{im}(du) \;,\\
\Phi_{2n}(\beta_0) & = & - {1 \over \sqrt n} \sum_{i=1}^n \sum_m
\int_0^{\tau_0} 
[Z_{im}(u)  - { s^{(1)} \over s^{(0)}}(u, \beta_0, X_{im}) ] 
\quad \times \\
& & \quad \quad \quad \times \quad 
({S^{(0)}_{-i} -s^{(0)} \over s^{(0)}} )(u,\beta_0, X_{im}) 
N_{im}(du) \;,\\
\Phi_{3n}(\beta_0) & = & 
{1 \over \sqrt n} \sum_{i=1}^n \sum_m \int_0^{\tau_0}
({ S^{(1)}_{-i} - s^{(1)} \over s^{(0)}} ) (
{ S^{(0))}_{-i} - s^{(0)} \over s^{(0)}} ) (u, \beta_0,X_{im}) N_{im}(du) 
\;,\\ 
\Phi_{4n}(\beta_0) & = &  -
{1 \over \sqrt n} \sum_{i = 1}^n  
\sum_m \int_0^{\tau_0} {S^{(1)}_{-i} \over s} 
[({S^{(0)}_{-i} \over s^{(0)}} - 1 )^2 
{s^{(0)} \over S^{(0)}_{-i}} ](u,\beta_0,
X_{im}) N_{im}(du) \;, \\
\Sigma_{1n}(\beta) & = & {1 \over n} \sum_{i=1}^n \sum_m \int [{
S_{-i}^{(2)} - s^{(2)} \over s^{(0)}} - \psi_{-i} - \psi_{-i}^T](u, \beta,
X_{im}) N_{im}(du) \;,\\
\Sigma_{2n}(\beta) & = & - {1 \over n} \sum_{i=1}^n \sum_m \int [{
S_{-i}^{(1)} - s^{(1)} \over s^{(0)}}]^{\otimes 2}(u, \beta,
X_{im}) N_{im}(du) \;,\\
\Sigma_{3n}(\beta) & = & {1 \over n} \sum_{i=1}^n \sum_m \int
[\hat V_{-i} - \tilde V_{-i}] (u, \beta, X_{im}) N_{im}(du) \;,
\eeaa
where $\hat \psi_{-i} = (S_{-i}^{(1)} - s^{-1}) \otimes s^{(1)} /[s^{(0)}]^2$,
$\hat V_{-i} = S_{-i}^{(2)} /S_{-i}^{(0)} - 
(S_{-i}^{(1)} /S_{-i}^{(0)})^{\otimes 2}$ and
$\tilde V_{-i} = S_{-i}^{(2)} /s^{(0)} - 
(S_{-i}^{(1)} /s^{(0)})^{\otimes 2}$.

Under assumptions of Proposition 2.3, $\tilde \Sigma_{n}(\beta)$ 
is the negative
derivative of the score function $\tilde \Phi_n(\beta)$. Similarly, 
under assumptions of Proposition 2.5, we have
$\Phi_n(\beta_0) = \Sigma_{j=1}^4 \Phi_{jn}(\beta_0)$ and
$\Sigma_n(\beta) = \sum_{j=1}^3  
\Sigma_{nj}(\beta)$ is the negative
derivative of the score function $\Phi_n(\beta)$. 
The proof of both propositions
amounts to  application of the following lemma and results
of  Bickel et al ((1993), p. 517).

\blem 

\bitem
\item[{(i)}] Under assumptions of Proposition 2.3 we have
$\tilde \Phi_{0n}(\beta_0) \implies \N(0, \Sigma_2(\beta_0))$,
$\tilde \Sigma_n(\beta_0) \to_P \Sigma_1(\beta_0)$, $\tilde \Phi_{n}(\beta_0)
- \tilde \Phi_{0n}(\beta_0) \to_P 0$, and \newline
$\sup \{|\tilde \Sigma_n(\beta) - \tilde \Sigma(\beta_0)|: |\beta - \beta_0|
\leq \epsilon_n\} \to_P 0$.

\item[{(ii)}] Under assumptions of Proposition 2.5 we have
$\Phi_{0n}(\beta_0) \implies \N(0, \Sigma(\beta_0))$, \newline
$\Sigma_{0n}(\beta_0) \to_P
\Sigma(\beta_0)$, $\Phi_{1n}(\beta_0) - \Phi_{0n}(\beta_0) \to_P 0$, 
$\Phi_{kn}(\beta_0) \to_P 0$
for $k = 2,3,4$, $\Sigma_{kn}(\beta_0) \to_P 0$ for $k = 1,2,3$, and
$\sup \{| \Sigma_{kn}(\beta) -  \Sigma_{kn}(\beta_0)|: |\beta - \beta_0|
\leq \epsilon_n\} \to_P 0$ for $k = 0,1,2,3$.
\eitem
\elem

\Proof 
First note that under the assumed regularity conditions, 
asymptotic normality of the terms $\tilde \Phi_{0n}(\beta_0)$ and
$\Phi_{0n}(\beta_0)$ follows from CLT.

We  show that $\tilde \Phi_n(\beta_0) - 
\tilde \Phi_{0n}(\beta_0) \to_P 0$
and $\Phi_{1n}(\beta_0) -  \Phi_{0n}(\beta_0) \to_P 0$. For any
bounded function $\phi(u,x)$, let 
$G_{ij}^{\phi} = G^{\phi} (W_i, W_j)$ be given by
$$
G_{ij}^{\phi} =  \sum_m \int_0^{\tau_0} \phi(u,X_{im}) [Z_{im}(u)
S^{(0)}_j(u,\beta, X_{im}) -
S^{(1)}_j(u,\beta, X_{im})] N_{im}(du)
$$
Under assumptions of Proposition 2.3, we have
$\tilde \Phi_n(\beta_0) = \tilde \Phi_{0n}(\beta_0) + O_P(\sqrt n a) + 
\empU_{n,2}(\pi_2 G^{\phi}))$ for $\phi \equiv 1$. Similarly, under assumptions
of Proposition 2.5, we have $\Phi_{1n}(\beta_0) = 
\Phi_{0n}(\beta_0) + O_P(\sqrt n a^2) +
\empU_{n,2}(\pi_2 G^{\phi})$ for $\phi(u,x) = [s^{(0)}(u, \beta_0,x)]^{-1}$.
Thus it is enough to show that 
in both cases $\empE \empU_{n,2}(\pi_2 G^{\phi}) = O((na)^{-1/2})$.
Choose $\phi = [s^{(0)}]^{-1}$ for instance, and define
\beaa
& & \bar G_n(W_i,W_j))   = 
a^{-1}\sum_{p=0}^1 \sum_{m} \int |Z_{im}(u)|^{p} 
\bar f_{1-pjn}(u, \beta_0,X_{im}) 
M_{im}(du) \\
&  + & a^{-1}
\sum_{p=0}^1 \int 
\sum_{m}  |Z_{im}|^p(u) Y_{im}(u) e^{\beta_0 Z_{im}(u)}
\bar f_{1-p,jn}( u, \beta_0,X_{im})
\alpha(u, X_{im})du \\
& & 
 a^{-1}\sum_m  \int_0^{\tau_0}
[|Z_{im}(u)|\bar f_{0nj}(u,\beta_0, X_{im}) + 
\bar f_{1,nj}(u,\beta_0, X_{im})] 
N_{im}(du)  \;.
\eeaa
We have $\empE \empU_{n,2}(\pi_2 G^{\phi})) = 
O(n^{-1/2}( \empE\bar G^2_n(W_1,W_2))^{1/2}) = O((na)^{-1/2})$ because, by
Lemma 2.1, the expectation $\empE \bar G_n^2(W_1,W_2)$ is bounded by  
\beaa
& & {4 \over a^2}  \sum_{p=0}^1 \int_0^{\tau_0} \int_0^{\tau}
\sigma_{p,p}(u,u,x) \empE [\bar f_{1-p,jn}(u,\beta_0,x)]^2
\alpha(u,x) du  dx +  \\  
& & \int_0^{\tau_0}\int_0^{\tau_0} \int_0^{\tau}
 \sigma_{p,p}(u_1,u_2,x)
\empE [\prod_{l=1}^2 \bar f_{1-p,jn}(u_l, \beta_0,x)] 
\prod_{l=1}^2 \alpha(u_l, x) du_1 du_2 dx + \\
& &   
\int_0^{\tau_0} \int_0^{\tau_0} \int_0^{\tau} \int_0^{\tau}
\rho_{p;p}(\underline u, \underline x ) 
\empE [\prod_{l=1}^2 \bar f_{1-p,jn}(u_l, \beta_0,x_l)] 
\prod_{l=1}^2 \alpha(u_l, x_l) du_l dx_l 
\eeaa
Here in the last line
$\underline u = (u_1, u_2)$ and $\underline x = (x_1, x_2)$.
By Lemma 6.8, the bound
is of order $O(a^{-1})$. 
It follows now that $\Phi_{1n}(\beta_0) - \Phi_{0n}(\beta_0) \to_P 0$. 

The same argument  applied to the function $\phi(u,x) \equiv 1$ shows 
that $\tilde \Phi_n(\beta_0) - \tilde \Phi_{0n} \to_P 0$.  Changing
the risk processes $S_j^{(k)}$ $-S_j^{(k+1)}$, $k=0,1$, in the definition
of $G^{\phi}(W_i,W_j)$, we also obtain 
\beaa
& & \tilde  \Sigma_n(\beta_0) - \Sigma_1(\beta_0) = \\
& = & n^{-1} \sum_m \int_0^{\tau_0}
Z_{im}(u)  s^{(1)}(u, \beta_0, X_{im})
-  s^{(2)}(u, \beta_0, X_{im}) M_{im}(du) + o_p(1)
\eeaa
The Strong Law of Large Numbers implies that
 $\tilde  \Sigma(\beta_0) \to_P \Sigma_1(\beta_0)$.
Components of the matrix $\tilde \Sigma(\beta)$ are Lipschitz
continuous in $\beta$, and it is easy to verify that  
$|\tilde \Sigma_n(\beta) - \tilde \Sigma_n(\beta')| \leq 
|\beta - \beta'| \empU_{n,2}(G_{2n})$ where 
$G_{2n}$ is a kernel degree 2 satisfying $\empE |\empU_{n,2}(G_{2n})| = 
O(1)$. This completes the proof of the first part of the 
proposition.

Further, the terms $\Phi_{2n}(\beta_0)$ and $\Sigma_{1n}(\beta_0)$
are U-statistics of degree 2. Using similar algebra as in the case
of the difference $\Phi_{1n} - \Phi_{0n}$, we can show that they
converge to 0 in probability.

Next define
$$
H_{1n}  =  {1 \over n} \sum_{i=1}^n \sum_m \int_0^{\tau_0} \phi(u, X_{im}) \prod_{k=1}^2
({S_{-i}^{(p)} -
s^{(p)} \over s^{(0)}} 
{S_{-i}^{(q)} -
s^{(q)} \over s^{(0)}} 
) (u, \beta_0,X_{im}) N_{im}(du)
$$
where $\phi(u,x)$ is a bounded function and $p,q  = 0$ or  $1$.
We have $\sqrt n H_{1n} = O_p(\sqrt n a^2) + O(1) 
[\sqrt n a^{-2} \empU_{n,3}(H) +
\sqrt n (na^2)^{-1} \empU_{n,2}(\bar H)]$, where
$$
H(W_i,W_j,W_k)  =  \sum_m \int_0^{\tau_0} \phi(u, X_{im})[f_{pjn} - \empE f_{pjn}]
[f_{qkn} - \empE f_{qkn}]
(u,\beta_0,X_{im}) N_{im}(du) 
$$
and $\bar H(W_i,W_j) = H(W_i,W_j,W_j)$.
We have $\empE H(W_1, W_2, W_3) = 0$. Lemmas 5.7, 6.8 and 6.9 imply
that  
$\empE \sqrt na^{-2} |\empU_{n,3}(\pi_3 H)| = O((na)^{-1})$ and
$O((na)^{-1/2}) = $ \newline 
$\empE \sqrt n (na^{2})^{-1} |\empU_{n,2}(\pi_2 \empE_{\{23\}} H)|$
while the remaining projections are 0. Further, 
$\sqrt n (na^2)^{-1} \empE \empU_{n,2} (\bar H) = O((na^2)^{-1/2}) =
\sqrt n (na^2)^{-1} \empU_{n,2}(|\bar H|)$ so that the condition $na^2 \uparrow \infty$ implies asymptotic negligibility of the third term of 
$\sqrt n H_{1n}$.

The choice of $\phi \equiv 1$, $p = 1, q = 0$ implies that 
if $na^4 \downarrow 0$ and $na^2 \uparrow \infty$
then $\Phi_{3n}(\beta_0) \to_P 0$. 
The choice of $\phi \equiv 1$ and $p = q =1$ implies 
$\Sigma_{2n}(\beta_0) \to_P 0$.

To handle the term $\Phi_{4n}(\beta_0)$ define
$$
H_{2n}  =   {1 \over n} \sum_{i=1}^n
\sum_m \int_0^{\tau_0} ({\bar f_{1n} \over s^{(0)}} 
 ({S_{-i}^{(0)} -
s^{(0)} \over s^{(0)}} )^2 ) (u, \beta_0,X_{im}) N_{im}(du) \;.
$$
Using $(x+y)^2 \leq 2x^2 + 2y^2$, we have
$\sqrt n |H_{2n}|  \leq  2 O_p(\sqrt n a^4) + 2\sqrt n H_{2n;1} +
 2\sqrt n H_{2n;2}$, 
where $H_{2n;1}$ corresponds the sum $H_{1n}$ applied with function 
$\phi = \empE \bar f_{1ni}/ s^{(0)}$, and $H_{2n;2}$ is a 
V statistics of degree 4: $H_{2n;2} = O(1)[a^{-3} \empU_{n,4}(h) + 
(a^3 n)^{-1}
\empU_{n,3}(\bar h) + 2(na^3)^{-1} \empU_{n,3}(h') + 
(n^2 a^3)^{-1} \empU_{n,2}(h'')]$,
where 
$$
h(W_i,W_j,W_k,W_l)  = 
 \sum_m \int_0^{\tau_0} [\bar f_{1jn} - \empE \bar f_{1jn}]
\prod_{p=k,\l}[f_{0pn} - \empE  f_{0pn}] (u,\beta_0,X_{im})
N_{im}(du) 
$$
and $\bar h(W_i,W_j,W_k) = h(W_i,W_j,W_k, W_k)$, 
$ h'(W_i,W_j,W_k) = h(W_i,W_j,W_j, W_k)$, $h''(W_i,W_j) = h(W_i,W_j,W_j,W_j)$.
We have
$ \empE 
|{\sqrt n (n^2a^3)^{-1}} \empU_{n,2} (h'')| \leq \newline
 {\sqrt n(na^3)^{-1}}  \empE \empU_{n,2}|h''|$ which is bounded by 
$$
 {\sqrt n \over 
 n^2a^3}       \int_{\R}\empE |[
  f_{1jn} - \empE \bar f_{1jn}]
[f_{0jn} - \empE  f_{0jn}]^2|(u,\beta_0,x) 
s^{(0)}(u,\beta_0,x) \alpha(u,x)du dx 
$$
Under conditions D.2 (ii)
this bound is of order  $O(n^{-3/2} a^{-2})$ and
tends to 0 
if $na^2 \uparrow \infty$. 
A similar argument shows also that the second and third term of 
$\sqrt n H_{2n;2}$ have expectation tending to 0 when $na^2 \uparrow \infty$ and $na^4 \downarrow 0$.
The first term has expectation 0. By Lemmas 
5.7 and 6.9, we have
$\empE \sqrt n a^{-3}  |U_{4n}
 \pi_4 h | = O( (na)^{-3/2})$,
$\empE \sqrt n a^{-3}  |\empU_{n,3}(
\pi_3 \empE_{\{234\}} h) | = O((na)^{-1})
$,
while the remaining projections are 0.

Further,  for $\epsilon \in (0,1)$, define
$$
\Omega_n(\epsilon) = \{ {1 \over 1+\epsilon} \leq \min_i \inf_{(u,x) \in \R 
\atop \beta \in \B}
{s^{(0)} \over S^{(0)}_{-i}}(u,\beta,x) \leq  
\max_i \sup_{(u,x) \in \R \atop \beta \in \B} 
{s^{(0)} \over S^{(0)}_{-i}}(u,\beta,x) \leq  
{1 \over 1 - \epsilon} \}
$$
As in the proof of Proposition 3.4, the condition C implies
$P(\Omega_n(\epsilon)) \to 1$.
Also on the event $\Omega_n(\epsilon)$,
the term $\Phi_{4n}(\beta_0)$ satisfies
$ |\Phi_{4n}(\beta_0| \leq (1 - \epsilon)^{-1} \sqrt n H_{2n}$.
For any $\eta > 0$, we have
$P(|\Phi_{4n}(\beta_0)| > \eta) \leq 
P(\sqrt n H_{2n} > \eta, \Omega(\epsilon)) + 
P(\Omega_n^c(\epsilon)) 
 \leq  P(\sqrt n H_{2n} > (1-\epsilon)\eta) +
P(\Omega_n^c(\epsilon)) \to 0$.

Application of the condition C shows also that 
$\Sigma_{3n}(\beta_0) \to_P 0$. Finally, it is easy to verify that 
the matrices $\Sigma_{nk}, k = 0,1,2,3$, satisfy 
$|\Sigma_{nk}(\beta) - \Sigma_{nk}(\beta_0)| \leq |\beta - \beta'|O_P(1)$, 
which completes the proof of the lemma.

\par 

\noindent {\large\bf Acknowledgment}
We thank an anonymous reviewer and Editor Jane Ling Wang for 
comments. 
Research was supported by grants from the National Science Foundation and
National Cancer Institute. The data presented here were obtained from
the Statistical Center of the International Blood and Marrow Transplant
Registry. The analysis has not been reviewed or approved by the Advisory
Committee of the IBMTR. 

\par

\noindent{\large\bf References}

\begin{list}{0}{\setlength{\rightmargin}{\leftmargin}}

\item[
Andersen, P.K., Borgan, O., Gill, R.D. and Keiding, N.] (1993). 
{\it Statistical Models Based on Counting Processes}.  Springer
Verlag, New York. 
\item[
Andersen, P.K. and  Gill, R.D. ] (1982).  Cox's regression model
for counting processes: a large sample study. \it Ann. Stat. \bf 10
\rm 1100-1120.

\item[
Arjas, E. and Eerola, M.] (1993). On predictive causality in longitudinal
studies. \it J. Statist. Planning and Inference, \bf 34\rm,
361--384.
 
\item[
Beran, R. ] (1981).  Nonparametric regression with randomly
censored survival data. Tech. Report, University of
California, Berkeley.

\item[
Bickel, P.J., Klaassen, C., Ritov, Y. and Wellner, J.A. ] (1993). 
\it Efficient and Adaptive Estimation in transformation
models. \rm Johns Hopkins Univ. Press.

\item[
Cox, D. R.] (1973). The statistical analysis of dependencies
in point processes. \it Symposium on Point Processes \rm
(Lewis, P. A. W., Ed.). Wiley, New York.

\item[
Dabrowska, D. M., Sun, G.W. and Horowitz, M. M.] (1994). Cox
regression
in a Markov renewal model: an application to the analysis of bone
marrow
transplant data. \it J. Amer. Statist. Assoc. \bf 89\rm,
867--877.

\item[
Dabrowska, D.M.] (1997). Smoothed Cox regression. \it Ann. Statist.
\bf 25\rm,  1510-1540.

\item[
de la Pe\~na, V. and Gin\'e, E] (1999). \it Decoupling: From Dependence
to Independence. \rm Springer Verlag.

\item[
Einmahl, U. and Mason, D.] (2000). An empirical process approach to
uniform consistency of kernel-type function estimators.  \it J. Theor.
\newline Probab. \bf 13\rm,  1-37.

\item[
Gill, R. D.] (1980). Nonparametric estimation based on censored
observations of a Markov renewal process. \it Z.
Wahrscheinlichkeitstheorie verv. Gebiete, \bf 53\rm, 97-116.

\item[
Gin\'e, E. and Guillou, A.] (1999) Laws of iterated logarithm
for censored data. \it Ann. Probab. \bf 27\rm, 2042-2067.

\item[
Jacod, J.] (1975). Multivariate point processes: predictable projection,
Radon-Nikodym derivatives, representation of martingales. \it Z.
\newline Wahrscheinlichkeitstheorie verv. Gebiete 
 \bf 31\rm, 235-254.

\item[
Klein, J.P., Keiding, N. Copelan, E.A.] (1993).
Plotting summary predictions in multistate survival models: probabilities
of relapse and death in remission for bone marrow transplant patients.
\it Statistics in \newline Medicine \bf 12\rm, 2315-2332.

\item[
Lagakos, S.W., Sommer, C.J. and Zelen, M.] (1978) Semi-Markov models
for partially censored data. \it Biometrika \bf 65\rm, 311-317.

\item[
M\"uller, H.G.  and Wang, J.L.] (1994) Hazard rate estimation under random 
censoring with varying kernels and bandwidths. \it Biometrics\bf
50\rm, 61-76.

\item[
Nielsen, J. P., Linton, O.B. and Bickel, P.J.] (1998) On a semi-parametric 
survival model with flexible covariate effect. \it Ann. Statist.
\bf 26\rm, 215-241.

\item[
Nolan, D. and Pollard, D.] (1987). U-processes: rates of  convergence. \newline
\it Ann. Statist. \bf 15\rm, 780-799.

\item[
Oakes, D. and Cui, L.] (1994). On semi-parametric inference
for modulated renewal processes. \it Biometrika \bf 81\rm, 83-90.

\item[
Pakes, A. and Pollard, D.] (1989) Simulation and the asymptotics
of the optimization estimators. \it Econometrica \bf 57\rm, 1027-1057.

\item[ 
Phelan, M.J.] (1990). Bayes estimation from a Markov renewal process.
\it Ann. Statist. \bf 18\rm, 603-616.

\item[
Pons, 0. and Visser, M.] A non-stationary Cox model. \it Scand. J. Statist.
\bf 27, \rm 619--641.
 
\item[
Sasieni, P.] (1992) Information bounds for conditional hazard ratio
in a nested family of regression models. \it J. Roy. Statist. Soc.
B \bf 54\rm, 617-635.

\item[
Van der Vaart, A. W. and Wellner, J. A.] (1996). \it Weak Convergence 
of \newline Empirical Processes with Applications. \rm Springer Verlag.

\item[
Voelkel, J. G. and Crowley, J. J.] (1984). Nonparametric inference
for a class of semi-Markov processes with censored observations.
\it Ann. \newline Statist. \bf 12\rm,
142-160.

\end{list}

\vskip .65cm
\noindent
Dorota M. Dabrowska\\
Department of Biostatistics\\
University of California \\
Los Angeles, CA 90095-1772
\vskip 2pt
\noindent
E-mail: (dorota@brahms.ph.ucla.edu)
\vskip 2pt
\noindent
Wai Tung Ho \\
SPSS Inc., \\
233 South Wacker Drive, 11th Floor \\
Chicago, IL 60606
\vskip 2pt
\noindent
E-mail: (ho-wai-tung@alumni.cuhk.net)
\vskip .3cm

\newpage

%
%
%
%

\begin{center}
Table 3.1. Polynomial kernels of degree 2, 4 and 6
\end{center}
\begin{table}[htb]
\begin{tabular}{cccc}
\hline
$\mu=1$ & interior & $(3/4)(1-x^2)$ \\
        & left     & $6(p+x)(p+q)^{-4}[p^2-2pq+3q^2+2x(p-q)]$\\
        & right    & $6(q-x)(p+q)^{-4}[3p^2-2pq+q^2+2x(2p-q)]$\\
$\mu=2$ & interior & $(15/16)(1-x^2)^2$ \\
        & left     & $60(q-x)(p+x)^2(p+q)^{-6}[p^2-2pq+2q^2+(2p-3q)x]$\\
        & right    & $60(q-x)^2(p+x)(p+q)^{-6}[2p^2-2pq+q^2+(3p-2q)x]$\\
$\mu=3$ & interior & $(35/32)(1-x^2)^3$\\
        & left     & $140(q-x)^2(p+x)^3[3p^2-6pq+5q^2+2(3p-4q)x]$\\
        & right    & $140(q-x)^3(p+x)^2[5p^2-6pq+3q^2+2(4p-3q)x]$\\ 
\hline
\end{tabular}
\end{table}

\newpage

\begin{center}
Table 4.1. Regression estimates and standard errors of direct
transitions
\end{center}
\begin{table}[ht]

\centerline{
\begin{tabular}{cccccc} 
& & & & & \\
\hline
\hline
& {TX $\to$ AGVHD} & &{TX $\to$ CGVHD} 
& & {AGVHD $\to$ CGVHD} \\
\hline
\hline
& & & & &    \\
sex--match   & .08 (.05) & & .12 (.05) & &   \\
CSA & .46 (.08) & & .18 (.12) & &  \\
Trem & - .58  (.13) & &  -.42 (.15) & & -.28 (.24) \\
MTX	& .38 (.20) &  & -.40. (.31) &  &  \\
disease & .12 (.07) & & & & -.13 (.13) \\
\hline
\hline 
& {TX $\to$ relapse} & & {AGVHD $\to$ relapse} & & {CGVHD $ \to$ relapse}\\
\hline
\hline
& & & & &  \\
sex--match & -.11 (.10) & & & & .21 (.10) \\
CSA & & & -.52 (.31) & & \\
Trem & .20 (.12) & & -.75 (.59) & & .40 (.31) \\
MTX & .32 (.23) & & .67 (.56) & & \\
disease & .14 (.10) & & & & \\ 
\hline
\hline
& {TX $\to$ death} & & {AGVHD $\to$ death} & & {CGVHD $ \to$ death}\\
\hline
\hline
& & & & & \\
Trem & .23 (.15) & & .57 (.20) & & .48 (.25) \\
CSA &-.25 (.18) & & & & \\
MTX & -1.06 (0.58) & & & & -.82 (.71) \\
disease & & & .21 (.13) & & \\
prior AGVHD & & & & & .75 (.16) \\
\hline
\hline
\end{tabular} }
\end{table}

\noindent
The covariates are  binary 0-1 variables: 
Sex--match $ = 1$ if the donor is a female and the recipient
is a male. Disease $ = 1$ if the disease type is ALL;
Prior AGVHD $= 1$
if AGVHD occurs prior to  CGVHD.
The  GVHD prophylactic treatments are labeled as
cyclosporin (CSA =1), T cell removal (Trem = 1) and
methotraxate (MTX=1).

\newpage

Figure captions

Figure 4.1 Baseline cumulative hazards of transitions originating from
the transplant state versus age. The labels of states are
1 - transplant (TX), 2 - AGVHD, 3 - CGVHD, 4 - relapse and 5 - death.

Figure 4.2 Baseline cumulative hazards of transitions originating from
the AGVHD state versus age. The labels of states are
2 -AGVHD, 3 - CGVHD, 4 -
relapse and 5 -death.

Figure 4.3 Baseline cumulative hazards of transitions originating from
the CGVHD state versus age. The labels of states are 3 - CGVHD, 4 -
relapse and 5 -death.

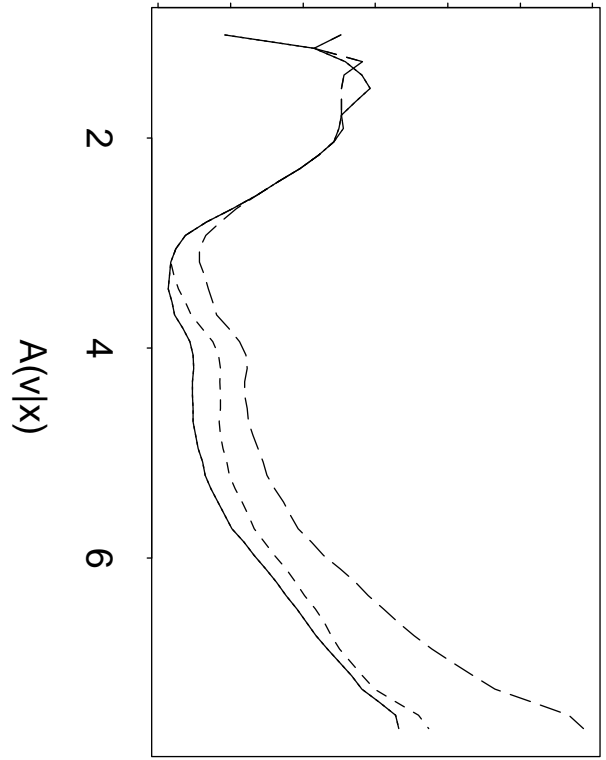
\end{document}